\documentclass[twoside,11pt,reqno]{amsart}
\usepackage{amsmath,amssymb,amscd,mathrsfs}

\makeatletter

\@addtoreset{equation}{section}

\def\thesubsection{\S\thesection-\alph{subsection}}
\def\Point#1{\addtocounter{subsection}{1}\vspace{2mm}
\noindent\thesubsection. {\bf #1.}\def\@currentlabel{\thesubsection}}

\newtheorem{Lemma}[equation]{Lemma}
\newtheorem{Theorem}[equation]{Theorem}
\newtheorem{Corollary}[equation]{Corollary}

\newtheorem{Remark}[equation]{Remark}

\newtheorem{Example}[equation]{Example}

\newtheorem{Conjecture}[equation]{Conjecture}

\newtheorem{Procedure}[equation]{Procedure}

\def\RR{{\mathtt{R}}}

\def\C{{\mathbb C}}
\def\Q{{\mathbb Q}}
\def\N{{\mathbb N}}
\def\Z{{\mathbb Z}}
\def\Mtype{\mathtt{M}}
\def\Qtype{\mathtt{Q}}
\def\0{{\bar 0}}
\def\1{{\bar 1}}
\def\pr{{\operatorname{pr}}}

\def\P{{p}}

\def\V{{\mathscr V}}
\def\T{{\mathscr T}}
\def\TT{{\widehat{\mathscr T}}}
\def\EE{{\widehat{\mathscr E}}}
\def\FF{{\widehat{\mathscr F}}}

\def\K{{\mathscr K}}
\def\KK{{\widehat{\mathscr K}}}
\def\E{{\mathscr E}}
\def\F{{\mathscr F}}
\def\U{{\mathscr U}}

\def\u{{\mathfrak u}}

\def\A{{\mathscr A}}

\def\cont{{\operatorname{cont}}}

\def\hom{{\operatorname{Hom}}}
\def\Ext{{\operatorname{Ext}}}
\def\End{{\operatorname{End}}}

\def\Lie{{\operatorname{Lie}}}

\def\ind{{\operatorname{ind}}}

\def\res{{\operatorname{res}}}

\def\soc{{\operatorname{soc}}}
\def\ch{{\operatorname{ch}\:}}
\def\wt{{\operatorname{wt}}}

\def\tr{{\operatorname{tr}\,}}
\def\Tr{{\operatorname{Tr}}}

\def\cont{{\operatorname{cont}}}
\def\underbar{\mathpalette\@underbar}
\def\@underbar#1#2{\settowidth{\@tempdimb}{$#1#2$}\@tempdimb=0.8\@tempdimb
                   \ooalign{$#1#2$\crcr%
                         \hfil\rule[-.5mm]{\@tempdimb}{.4pt}\hfil}}

\def\eps{{\varepsilon}}
\def\phi{{\varphi}}

\def\la{{\lambda}}

\def\ev{{\operatorname{ev}}}

\newdimen\hoogte    \hoogte=12pt    
\newdimen\breedte   \breedte=14pt  
\newdimen\dikte     \dikte=0.5pt 

\newenvironment{Young}{\begingroup
       \def\vr{\vrule height0.89\hoogte width\dikte depth 0.2\hoogte}
       \def\fbox##1{\vbox{\offinterlineskip
                    \hrule height\dikte
                    \hbox to \breedte{\vr\hfill##1\hfill\vr}
                    \hrule height\dikte}}
       \vbox\bgroup \offinterlineskip \tabskip=-\dikte \lineskip=-\dikte
            \halign\bgroup &\fbox{##\unskip}\unskip  \crcr }
       {\egroup\egroup\endgroup}
\def\diagram#1{\relax\ifmmode\vcenter{\,\begin{Young}#1\end{Young}\,}\else%
              $\vcenter{\,\begin{Young}#1\end{Young}\,}$\fi}

\begin{document}
\title[Character formulae for $\mathfrak q(n)$]{\boldmath
Kazhdan-Lusztig polynomials and character formulae for the Lie superalgebra $\mathfrak{q}(n)$}
\author{\sc Jonathan Brundan}
\address
{Department of Mathematics\\ University of Oregon\\
Eugene\\ OR~97403, USA}
\email{brundan@darkwing.uoregon.edu}
\thanks{Work
partially supported by the NSF (grant no. DMS-0139019).}
\maketitle

\vspace{-2mm}
\section{Introduction}

The problem of computing the characters of the finite dimensional
irreducible representations for the classical Lie superalgebras
was posed originally by Kac in 1977 \cite{Kac, Kac2, Kac3}.
For the family $\mathfrak{q}(n)$,
various special cases 
were treated by Sergeev \cite{sergeev} (polynomial representations)
and Penkov \cite{penkov, Penk} (typical then generic representations),
culminating in the complete solution of Kac' problem for 
$\mathfrak{q}(n)$
in the work of Penkov and Serganova \cite{PS1, PS2} in 1996.
In this article, we will explain a different approach 
to computing the characters
of the irreducible ``integrable'' representations of 
$\mathfrak{q}(n)$, i.e. the representations that lift to 
the supergroup $Q(n)$. 

The strategy followed runs parallel to our recent work 
\cite{KL} on representations
of $GL(m|n)$, and is in keeping with the Lascoux-Leclerc-Thibon
philosophy \cite{LLT}.
We first study the canonical basis of the representation
$$
\F^n := {\bigwedge}^n \V,
$$
where $\V$ denotes the natural representation of 
the quantized enveloping algebra $U_q(\mathfrak{b}_{\infty})$.
This provides a natural Lie theoretic framework for the combinatorics
associated to the representation theory of $Q(n)$.
The idea that $\mathfrak{b}_{\infty}$ should be relevant here
is already apparent from \cite{LT, BK}.
Our main theorem shows that the transition matrix between the canonical basis
and the natural monomial basis of $\F^n$ at $q = 1$
is transpose to the transition matrix 
between the bases for the Grothendieck
group of finite dimensional representations of $Q(n)$
given by certain Euler characteristics and by 
the irreducible representations.

In order to define the canonical basis of $\F^n$,
we must also consider the tensor space 
$$
\T^n := {\bigotimes}^n \V.
$$
Work of Lusztig \cite[ch. 27]{Lubook} 
shows how to construct a canonical basis for $\T^n$.
We then pass from there to the space
$\F^n$, which we realize as a quotient of
$\T^n$ following 
Jing, Misra and Okado \cite{JMO}.
The entries of the transition matrix between the canonical basis and the natural monomial basis of $\T^n$ should be viewed
as the combinatorial analogues for the Lie superalgebra $\mathfrak{q}(n)$ 
of the Kazhdan-Lusztig polynomials of \cite{KaL}.
We conjecture that these polynomials evaluated at $q = 1$
compute the composition multiplicities of 
the Verma modules in the analogue 
of category $\mathcal O$
for $\mathfrak{q}(n)$, see \ref{mc} for a precise
statement.

We now state our main result precisely.
Let $\Z^n_+$ denote the set of all tuples
$\la = (\la_1,\dots,\la_n) \in \Z^n$
such that $\la_1 \geq \dots\geq \la_n$ and moreover
$\la_r = \la_{r+1}$ implies $\la_r = 0$ for each $r = 1,\dots,n-1$.
For $\lambda \in \Z_+^n$,
let $z(\lambda)$ denote the number of $\lambda_r \:(r=1,\dots,n)$
that equal zero.
Also let $\delta_r$ denote the $n$-tuple
with $r$th entry equal to $1$ and all other entries equal to zero.
Given $\lambda \in \Z^n_+$, there is an irreducible 
representation $L(\lambda)$ of $Q(n)$
of highest weight $\lambda$, unique up to isomorphism.
Let $L_\lambda$ denote the character of $L(\lambda)$, giving us a 
canonical basis
$\{L_\lambda\}_{\lambda \in \Z^n_+}$ for the character ring of $Q(n)$,
see \ref{reps}.
There is another basis denoted $\{E_\mu\}_{\mu \in \Z^n_+}$
which arises naturally from certain Euler characteristics,
see \ref{ecsec}.
We can write
$$
E_\mu = \sum_{\la \in \Z^n_+} d_{\mu,\la} L_\la
$$
for coefficients $d_{\mu,\la} \in \Z$, where $d_{\mu,\mu} = 1$
and $d_{\mu,\la} = 0$ for $\la\not\leq\mu$ in the dominance ordering.
The $E_\mu$'s are explicitly known: they are multiples
of the symmetric functions known as {\em Schur's $P$-functions}.
So the problem of computing $L_\lambda$ for each $\la \in \Z^n_+$
is equivalent to determining  the decomposition numbers
$d_{\mu,\la}$ for each $\la,\mu\in\Z^n_+$.

\vspace{1mm}
\noindent{\bf Main Theorem.}
{\em Let $\lambda = (\lambda_1,\dots,\lambda_n) \in \Z^n_+$.
Choose $p$ to be maximal such that there exist
$1\leq r_1 < \dots < r_p < s_p < \dots < s_1 \leq n$ 
with $\lambda_{r_q} + \lambda_{s_q} = 0$ for all
$q = 1,\dots,p$.
Let $I_0 =\{|\la_1|,\dots,|\la_n|\}$.
For $q = 1,\dots,p$,
define $I_q$ and $k_q$
inductively according to the following rules:
\begin{itemize}
\item[(1)]
if $\la_{r_q} > 0$, let $k_q$ be the smallest
positive integer with $\la_{r_q} + k_q \notin I_{q-1}$, and 
set $I_q = I_{q-1} \cup \{\la_{r_q} + k_q\}$;
\item[(2)]
if $\la_{r_q} = 0$, let $k_q$ and $k_q'$ be the smallest positive
integers with $k_q, k_q' \notin I_{q-1}$,
$k_q < k_q'$ if $z(\la)$ is even and
$k_q > k_q'$ if $z(\la)$ is odd, and set
$I_q = I_{q-1} \cup \{k_q,k_q'\}$.
\end{itemize}
Finally, for each $\theta = (\theta_1,\dots,\theta_p) \in \{0,1\}^p$,
let $\RR_{\theta}(\lambda)$ denote the unique element of $\Z^n_+$ that is 
conjugate to
$\la + \sum_{q=1}^p \theta_q k_q (\delta_{r_q} - \delta_{s_q})$.
Then, 
$$
d_{\mu,\la} = 
\left\{
\begin{array}{ll}
2^{(z(\lambda)-z(\mu))/2}&\hbox{if $\mu = \RR_{\theta}(\lambda)$ for some
$\theta = (\theta_1,\dots,\theta_p)\in\{0,1\}^p$},\\
0&\hbox{otherwise.}
\end{array}
\right.
$$
}
\vspace{1mm}

The remainder of the article is organized as follows.
In $\S$\ref{ta}, we introduce the quantum group of 
type $\mathfrak{b}_{\infty}$, 
and construct the canonical basis of the tensor space $\T^n$.
In $\S$\ref{ea}, we pass from $\T^n$ to the quotient 
$\F^n$, and study its canonical basis.
This time it turns out to be quite easy to compute explicitly.
Finally in $\S$\ref{cf}, we prove the character formulae.
Note there is one difficult place in our proof when we
need to appeal to the existence of certain homomorphisms between
Verma modules, see Lemma~\ref{ps}. 
For this we appeal to the earlier work of Penkov
and Serganova \cite[Proposition 2.1]{PS1}, which in turn relies upon
a special case 
of Penkov's generic character formula \cite[Corollary 2.2]{Penk}.
It would be nice to find an independent proof of this fact.

\vspace{1mm}
\noindent{\em Notation.}
Generally speaking, indices $i,j,k$ will run over
the natural numbers $\N = \{0,1,2,\dots\}$, indices
$a,b,c$ will run over all of $\Z$, and indices $r,s,t$
will run over the set $\{1,2,\dots,n\}$ where $n$ is a fixed
positive integer.

\section{Tensor algebra}\label{ta}

\Point{\boldmath Quantum group of type $\mathfrak{b}_{\infty}$}\label{qg}
Let $P$ be the free abelian group on basis
$\eps_1,\eps_2,\eps_3,\dots$, equipped with a symmetric bilinear form
$(.,.)$ defined by $(\eps_i, \eps_j) = 2\delta_{i,j}$
for all $i, j \geq 1$.
Inside $P$, we have the root system
$\{\pm\eps_i, \pm \eps_i \pm \eps_j\:|\:1 \leq i <  j\}$
of type $\mathfrak{b}_{\infty}$.
We use the following labeling for the Dynkin diagram:
$$
{\begin{picture}(100, 15)%
\put(6,5){\circle{4}}%
\put(12,2.3){$<$}%
\put(25, 5){\circle{4}}%
\put(44, 5){\circle{4}}%
\put(8, 4){\line(1, 0){15.5}}%
\put(8, 6){\line(1, 0){15.5}}%
\put(27, 5){\line(1, 0){15}}%
\put(46, 5){\line(1, 0){1}}%
\put(49, 5){\line(1, 0){1}}%
\put(52, 5){\line(1, 0){1}}%
\put(55, 5){\line(1, 0){1}}%
\put(58, 5){\line(1, 0){1}}%
\put(61, 5){\line(1, 0){1}}%
\put(64, 5){\line(1, 0){1}}%
\put(67, 5){\line(1, 0){1}}%
\put(70, 5){\line(1, 0){1}}%
\put(73, 5){\line(1, 0){1}}%
\put(6, 11){\makebox(0, 0)[b]{$_{0}$}}%
\put(25, 11){\makebox(0, 0)[b]{$_{1}$}}%
\put(44, 11){\makebox(0, 0)[b]{$_{{2}}$}}%
\end{picture}}
$$
We take the simple roots $\alpha_0,\alpha_1,\dots$ defined from
$$
\alpha_0 = -\eps_1,\qquad
\alpha_i = \eps_i - \eps_{i+1}\quad(i \geq 1).
$$
This choice induces a dominance ordering $\leq$ on $P$:
$\beta \leq \gamma$ if $\gamma - \beta$ is a sum of simple roots.
The Cartan matrix $A = (a_{i,j})_{i,j \geq 0}$ is defined by
$a_{i,j} = 2 (\alpha_i,\alpha_j) / (\alpha_i,\alpha_i)$.

We will work over the ground field $\Q(q)$, where $q$ is an indeterminate.
Let $q_i = q^{(\alpha_i,\alpha_i) / 2}$, i.e.
$q_0 = q$ and $q_i = q^2$ for $i > 0$. Define the quantum integer
$$
[n]_i = \frac{q_i^n - q_i^{-n}}{q_i-q_i^{-1}}
$$
and the associated quantum factorial
$[n]_i^! = [n]_i [n-1]_i \dots [1]_i$.
There is a field automorphism $-:\Q(q) \rightarrow \Q(q)$
with $\overline{q} = q^{-1}$. We will call an additive map
$f:V \rightarrow W$ between $\Q(q)$-vector spaces
{\em antilinear} if $f(c v) = \overline{c} f(v)$
for all $c \in \Q(q), v \in V$.

The quantum group $\U = U_q(\mathfrak{b}_\infty)$ 
is the $\Q(q)$-algebra generated by elements
$E_i, F_i, K_i\:(i \geq 0)$ subject to the relations
\begin{align*}
K_iK_i^{-1} = K_i^{-1}K_i = 1,
&\qquad K_iK_j = K_jK_i,\\
K_iE_jK_i^{-1} = q^{(\alpha_i, \alpha_j)} E_j,
&\qquad K_iF_jK_i^{-1} = q^{-(\alpha_i, \alpha_j)} 
F_j,\\
E_iF_j - F_jE_i &= \delta_{i, j} \frac{K_{i} - K_{i}^{-1}}{q_i - q_i^{-1}},\\
\sum_{k=0}^{1-a_{i,j}} (-1)^k E_i^{(k)} E_j E_i^{(1-a_{i,j}-k)}
&=
\sum_{k=0}^{1-a_{i,j}} (-1)^k F_i^{(k)} F_j F_i^{(1-a_{i,j}-k)}
=0\qquad(i \neq j),
\end{align*}
where $E_i^{(r)} :=E_i^r / [r]_i^!$ and
$F_i^{(r)} :=F_i^r / [r]_i^!$.
Also let
$$
\left[
\begin{array}{c}
K_i\\r
\end{array}
\right]
=
\prod_{s=1}^r \frac{K_i q_i^{1-s} - K_i^{-1} q_i^{s-1}}{q_i^s - q_i^{-s}}.
$$
for each $r \geq 1$.

We regard $\U$ as a Hopf algebra with comultiplication
$\Delta:\U \rightarrow \U \otimes \U$ defined on generators by
\begin{align*}
\Delta(E_i) &= 1 \otimes E_i + E_i \otimes K_{i}^{-1},\\
\Delta(F_i) &= K_{i} \otimes F_i + F_i \otimes 1,\\
\Delta(K_i) &= K_i \otimes K_i.
\end{align*}
This is the comultiplication used by Kashiwara \cite{Ka}, 
which is different from the one in Lusztig's book \cite{Lubook}.

Let us introduce various (anti)automorphisms of $\U$.
First, we have the bar involution $-:\U \rightarrow \U$,
the unique antilinear algebra automorphism such that
\begin{align}
\overline{E_i} &=E_i,
&\overline{F_i} &= F_i, 
&\overline{K_i} &= K_i^{-1}.\\\intertext{We will also need the linear algebra
antiautomorphisms $\sigma,\tau:\U \rightarrow \U$ 
and the linear algebra automorphism
$\omega:\U\rightarrow\U$ defined by}
\sigma(E_i) &= E_i,
&\sigma(F_i) &= F_i,
&\sigma(K_i) &= K_{i}^{-1},\\
\tau(E_i) &= q_i F_i K_i^{-1},
&\tau(F_i) &= q_i^{-1} K_i E_i,
&\tau(K_i) &= K_{i},\\
\omega(E_i) &= F_i,
&\omega(F_i) &= E_i,
&\omega(K_i) &= K_{i}^{-1}.
\end{align}

\begin{Lemma}\label{caut}
The maps $\tau$ and $-\circ\sigma$ are coalgebra automorphisms,
i.e. we have that 
$\Delta(\phi(x)) = (\phi \otimes \phi)(\Delta(x))$
for all $x \in \U$ and either $\phi = \tau$ or $\phi = -\circ\sigma$.
The map $\omega$ is a coalgebra antiautomorphism, i.e.
$\Delta(\omega(x)) = P((\omega \otimes \omega)(\Delta(x)))$
where $P$ is the twist $x \otimes y\mapsto y \otimes x$.
\end{Lemma}

We will occasionally need Lusztig's $\Z[q,q^{-1}]$-form
$\U_{\Z[q,q^{-1}]}$ for $\U$. We recall from
\cite[$\S$6]{LuUnity} that this is the 
$\Z[q,q^{-1}]$-subalgebra of $\U$ generated by the elements
$E_i^{(r)}, F_i^{(r)}$,  $K_i^{\pm 1}$ 
and $\left[\begin{array}{c}K_i\\r\end{array}\right]$ 
for all $i \geq 0$, $r \geq 1$.
It inherits from $\U$ the structure of a Hopf algebra over $\Z[q,q^{-1}]$.

\Point{\boldmath Tensor space}\label{ts}
The natural representation of $\U$ is the $\Q(q)$-vector space $\V$
on basis $\{v_a\}_{a \in \Z}$ with action defined by
\begin{align*}
E_0 v_a &= \delta_{a,0} (q+q^{-1}) v_{-1} + \delta_{a,1} v_0,\quad
&E_i v_a &= \delta_{a,i+1} v_{i} + \delta_{a, -i} v_{-i-1},\\
F_0 v_a &= \delta_{a,0} (q+q^{-1}) v_{1} + \delta_{a,-1} v_0,\quad
&F_i v_a &= \delta_{a,i} v_{i+1} + \delta_{a,-i-1} v_{-i},\\
K_0 v_a &= q^{2\delta_{a,-1} - 2\delta_{a,1}} v_a,\quad
&K_i v_a &= q^{2\delta_{a,i}- 2\delta_{a,i+1} + 2\delta_{a,-i-1} - 
2\delta_{a,-i}} v_a,
\end{align*}
for all $a \in \Z$ and $i \geq 1$,
see for example \cite[$\S$5A.2]{Jantzen2}.
Let $\T = \bigoplus_{n \geq 0} \T^n$ be the tensor algebra of $\V$,
so  $\T^n = \V \otimes\dots\otimes\V$ ($n$ times)
viewed as a $\U$-module in the natural way.

Let $\Z^n$ denote the set of all $n$-tuples
$\lambda = (\lambda_1,\dots,\lambda_n)$ of integers.
The symmetric group $S_n$ acts on $\Z^n$ by the rule 
$w  \lambda= (\lambda_{w^{-1} 1}, \dots, \lambda_{w^{-1} n})$.
We will always denote the longest element of $S_n$ by $w_0$, 
so $w_0 \lambda = (\lambda_n,\dots,\lambda_1)$.
Given $\lambda \in \Z^n$, let
\begin{equation}\label{ndef}
N_\lambda = v_{\lambda_1} \otimes \dots \otimes v_{\lambda_n} \in \T^n.
\end{equation}
The vectors $\{N_\lambda\}_{\lambda \in \Z^n}$ obviously give
a basis for $\T^n$.
For $\lambda \in \Z^n$, let $z(\lambda)$ denote the number of
$\lambda_r\:(r=1,\dots,n)$ that equal zero.
We get another basis $\{M_\lambda\}_{\lambda \in \Z^n}$ 
for $\T^n$ by defining
\begin{equation}\label{mdef}
M_\lambda = \left({q+q^{-1}}\right)^{-z(\lambda)} N_\lambda
\end{equation}
for each $\lambda \in \Z^n$.

Let $(.,.)$ be the symmetric bilinear form on $\T^n$
such that 
\begin{align}
(M_\lambda,N_\mu) &= q^{z(\lambda)} \delta_{\lambda,\mu}\\\intertext{for all $\lambda,\mu \in \Z^n$.
Also define an antilinear automorphism $\sigma:\T^n\rightarrow\T^n$
and a linear automorphism $\omega:\T^n\rightarrow\T^n$ 
by}
\sigma(N_\lambda) &= q^{-z(\lambda)}N_{-\lambda},\\
\omega(N_\lambda) &= N_{-w_0\lambda}
\end{align}
for all $\lambda \in \Z^n$.
The following lemma is checked
by reducing to the case $n = 1$ using Lemma~\ref{caut}.

\begin{Lemma}\label{tricks} The following hold for all $x \in \U$ and $u,v \in \T^n$:
\begin{itemize}
\item[(i)]
$\sigma(xv) = \tau(\overline{\sigma(x)})\sigma(v)$;
\item[(ii)]
$(x u, v)= (u, \tau(x) v)$;
\item[(iii)] $\omega(xu) = \omega(x) \omega(u)$.

\end{itemize}
\end{Lemma}

\Point{Bruhat ordering}\label{bo}
Recall that a vector $v$ in a $\U$-module is said to be of 
{\em weight} $\gamma \in P$ if
$K_i v = q^{(\alpha_i, \gamma)} v$ for each $i \geq 0$.
The weight of the basis vector $v_i$ of $\V$
is $\eps_i$, where we write
$\eps_0 = 0$ and $\eps_{-i} = - \eps_i$.
Hence the vectors $N_\lambda$ and $M_\lambda$ have weight
\begin{equation}\label{wtdef}
\wt(\lambda) := \eps_{\lambda_1} + \dots + \eps_{\lambda_n} \in P.
\end{equation}
More generally, for each $r = 1,\dots,n$, let
\begin{equation}
\label{wtdef2}
\wt_r(\lambda) := \eps_{\lambda_r} + \dots + \eps_{\lambda_n} \in P.
\end{equation}
The {\em Bruhat ordering} $\preceq$ on the set $\Z^n$ 
is defined by $\lambda \preceq \mu$ if
$\wt_r(\lambda) \leq \wt_r(\mu)$ for each $r = 1,\dots,n$ 
with equality for $r = 1$.
Here is an equivalent formulation of the definition.
Write $\lambda \downarrow \mu$ if 
one of the following holds:
\begin{itemize}
\item[(1)] 
for some $1 \leq r < s \leq n$ such that
$\lambda_r > \lambda_s$, we have that
$\mu_r = \lambda_s, \mu_s = \lambda_r$ and
$\mu_t = \lambda_t$ for all $t \neq r,s$;
\item[(2)] 
for some $1 \leq r < s \leq n$ such that $\lambda_r + \lambda_s = 0$,
we have that $\mu_r = \lambda_{r} - 1, \mu_s = \lambda_s + 1$
and $\mu_t = \lambda_t$ for all $t \neq r,s$.
\end{itemize}
The following lemma has been checked in many examples,
but I still have not found a proof in general. However, it 
is not used in an essential way in the remainder of the article.

\begin{Lemma}\label{bruhat}
$\lambda \succeq \mu$ if and only if there exist
$\nu_1,\dots,\nu_k \in \Z^n$ such that
$\lambda = \nu_1 \downarrow \nu_2 \downarrow \dots \downarrow \nu_k = \mu$.
\end{Lemma}

The {\em degree of atypicality} of $\lambda \in \Z^n$ is defined by
\begin{equation}\label{atyp}
\#\lambda := 
n - \frac{1}{2}\sum_{i \geq 1} |(\wt(\lambda), \varepsilon_i)|.
\end{equation}
We will say that $\lambda \in \Z^n$ is {\em typical}
if $\#\lambda \leq 1$. Equivalently, $\lambda$ is typical
if $\lambda_r +\lambda_s \neq 0$ for all
$1 \leq r < s \leq n$. 

\Point{\boldmath Bar involution}\label{bin}
We next define a bar involution on $\T^n$ (actually on some completion
$\TT^n$ of $\T^n$) compatible
with the bar involution on $\U$, following \cite[$\S$27.3]{Lubook}.
To begin with, let us recall the definition of Lusztig's
quasi-$R$-matrix $\Theta$, translated 
suitably since we are working with a different comultiplication.

Write $\U = \U^- \U^0\U^+ = \U^+ \U^0 \U^-$ for the usual
triangular decompositions of $\U$, so 
$\U^-$ is generated by all $F_i$, $\U^0$ is generated by all $K_i^{\pm 1}$
and $\U^+$ is generated by all $E_i$.
For a weight $\nu \geq 0$, 
let $\U^+_\nu$ (resp. $\U^-_\nu$) denote the subspace of $\U^+$
(resp. $\U^-$) spanned by all
monomials $E_{i_1} \dots E_{i_r}$ (resp $F_{i_1}\dots F_{i_r}$) with
$\alpha_{i_1}+\dots+\alpha_{i_r} = \nu$, so
$$
\U^+ = \bigoplus_{\nu \geq 0} \U^+_\nu,
\qquad
\U^- = \bigoplus_{\nu \geq 0} \U^-_{\nu}.
$$
Define $\tr \nu := \sum_{i \geq 0} a_i$ if 
$\nu = \sum_{i \geq 0} a_i 
\alpha_i$.
Let $(\U \otimes \U)^\wedge$ denote the 
completion of the vector space $\U \otimes \U$
with respect to the 
descending filtration $\left((\U \otimes \U)_{d}\right)_{d \in \N}$
where
$$
(\U \otimes \U)_{d} =
\sum_{\tr \nu \geq d}
\left(
\U^- \U^0 \U^+_\nu\otimes \U
 + 
\U \otimes \U^+ \U^0 \U^-_\nu\right).
$$
Exactly as in \cite[$\S$4.1.1]{Lubook}, 
we embed $\U \otimes \U$ into $(\U \otimes \U)^\wedge$ in the obvious way,
then extend the $\Q(q)$-algebra structure on $\U \otimes \U$ to
$(\U \otimes \U)^\wedge$ by continuity.
The bar involution on $\U \otimes \U$ is defined by
$\overline{x \otimes y} := \overline{x} \otimes \overline{y}$,
and also extends by continuity to $(\U \otimes \U)^\wedge$.
Finally, the antiautomorphism $\sigma \otimes \sigma: \U \otimes \U
\rightarrow \U \otimes \U$ 
and the automorphism $P \circ (\omega\otimes\omega):\U \otimes\U
\rightarrow \U\otimes \U$ (where $P$ is the twist $x \otimes y
\mapsto y \otimes x$)
extend to the completion too.

\begin{Lemma}\label{theta}
There is a unique family of 
elements $\Theta_\nu \in \U_\nu^+ \otimes \U_\nu^-$ such that $\Theta_0 = 1$ and 
$\Theta := 
\sum_\nu \Theta_\nu \in (\U \otimes \U)^\wedge$ satisfies
$\Delta({u}) \Theta =
\Theta \overline{\Delta(\overline{u})}$
for all $u \in \U$ 
(identity in $(\U \otimes \U)^\wedge$).
Moreover, each $\Theta_\nu$ belongs to
$\U_{\Z[q,q^{-1}]} \otimes_{_\Z[q,q^{-1}]} \U_{\Z[q,q^{-1}]}$.
\end{Lemma}

\begin{proof}
The first part of the lemma is \cite[Theorem 4.1.2(a)]{Lubook}.
The second part about integrality follows
using \cite[Corollary 24.1.6]{Lubook}.
Actually {\em loc. cit.}
applies only to finite type root systems
but one can pass from $\mathfrak{b}_n$ to $\mathfrak{b}_{\infty}$
by a limiting argument.
\end{proof}

\begin{Lemma}\label{tprops} The following equalities hold in $(\U \otimes \U)^\wedge$:
\begin{itemize}
\item[(i)]
$\Theta \overline{\Theta} = \overline{\Theta} \Theta = 1 \otimes 1$;
\item[(ii)]
$(\sigma \otimes \sigma)(\Theta) = \Theta$;
\item[(iii)]
$(\sigma \otimes \sigma)(\overline{\Theta}) = \overline{\Theta}$;
\item[(iv)]
$(P \circ (\omega \otimes \omega)) (\Theta) = \Theta$.
\end{itemize}
\end{Lemma}

\begin{proof}
(i) This is \cite[Corollary 4.1.3]{Lubook}.

(ii) Recall by Lemma~\ref{caut} that the map
$- \circ \sigma = \sigma \circ -$ is a coalgebra automorphism.
Using this one checks that
$$
\Delta(u) = (\sigma \otimes \sigma) (\overline{\Delta(\sigma(\overline{u}))}),\qquad
\overline{\Delta(\overline{u})} =
(\sigma \otimes \sigma) (\Delta(\sigma(u))).
$$
Now apply the antiautomorphism $\sigma \otimes \sigma$ to the 
equality
$\Delta(\sigma(u)) \Theta = \Theta \overline{\Delta(\sigma(\overline{u}))}$
to get that
$(\sigma \otimes \sigma)(\Theta) \overline{\Delta(\overline{u})}
= \Delta(u) (\sigma \otimes \sigma)(\Theta)$.
Hence, $(\sigma \otimes \sigma)(\Theta) = \Theta$ by the uniqueness in
Lemma~\ref{theta}.

(iii) Combine (i) and (ii).

(iv) Follows easily from the uniqueness and Lemma~\ref{caut}.
\end{proof}

Let $\TT^n$ denote the completion of the vector space $\T^n$ with respect
to the descending filtration 
$(\T^n_d)_{d \in \Z}$, where $\T^n_d$ is the subspace of 
$\T^n$ spanned by $\{N_\lambda\}$ for
$\la \in \Z^n$ with $\sum_{i=1}^n i \lambda_i \geq d$.
We embed $\T^n$ into $\TT^n$ in the natural way.
Note that $\TT^n$ contains all 
vectors of the form $N_\lambda + (*)$ where $(*)$ is an 
infinite linear combination of $N_\mu$'s with $\mu \prec \lambda$.
The action of $\U$ on $\T^n$ extends 
by continuity to $\TT^n$, as does the map 
$\omega:\T^n\rightarrow\T^n$.

Now we are ready to inductively define the
bar involution on $\TT^n$. 
It will turn out to satisfy the following properties:
\begin{itemize}
\item[(1)] $-:\TT^n\rightarrow\TT^n$ is a continuous, antilinear involution;
\item[(2)] $\overline{xv} = \overline{x}\,\overline{v}$
for all $x \in \U, v \in \TT^n$;
\item[(3)] 
$\overline{N_\lambda} \in
N_\lambda + \widehat\sum_{\mu \prec \lambda} \Z[q,q^{-1}] N_\mu$ 
and
$\overline{M_\lambda} \in
M_\lambda + \widehat\sum_{\mu \prec \lambda} \Z[q,q^{-1}] M_\mu$
for all $\lambda \in \Z^n$;
\item[(4)]
$\overline{\omega(v)} = \omega(\overline{v})$ for all $v \in \TT^n$.
\end{itemize}
If $n = 1$, the bar involution is defined by setting
$\overline{v_a} = v_a$ for each $a \in \Z$,
then extending by continuity.
The properties (1)--(4) in this case are easy to 
check directly.
Now suppose that $n > 1$, write $n = n_1+n_2$ for
some $n_1, n_2 \geq 1$, and assume 
we have already constructed bar involutions
on $\TT^{n_1}$ and $\TT^{n_2}$ satisfying properties (1)--(4).
Because of the way the completion $\TT^n$ is defined, 
multiplication by $\Theta$ gives a linear map
\begin{equation}
\Theta_{n_1,n_2}: \TT^{n_1} \otimes \TT^{n_2} \rightarrow \TT^{n}.
\end{equation}
Given $v \in \T^{n_1}, w \in \T^{n_2}$, 
let
$\overline{v \otimes w} := 
\Theta_{n_1, n_2} (\overline{v} \otimes \overline{w})$,
defining an antilinear map 
$-:\T^n \rightarrow \TT^n$.
The explicit form of $\Theta$ from Lemma~\ref{theta} combined with 
the inductive hypothesis implies that 
\begin{align*}
\overline{N_\lambda} &=
N_\lambda + \hbox{(a possibly infinite $\Q(q)$-linear combination of
$N_\mu$'s with $\mu \prec \lambda$)},\\
\overline{M_\lambda} &= 
M_\lambda + \hbox{(a possibly infinite $\Q(q)$-linear combination of
$M_\mu$'s with $\mu \prec \lambda$)}.
\end{align*}
Each $\Theta_\nu$ belongs
to $\U_{\Z[q,q^{-1}]} \otimes_{_{\Z[q,q^{-1}]}}
\U_{\Z[q,q^{-1}]}$ by Lemma~\ref{theta}, and $\U_{\Z[q,q^{-1}]}$ leaves
the $\Z[q,q^{-1}]$-lattices 
in $\TT^n$ generated by either the $N_\lambda$'s or the $M_\lambda$'s
invariant. Hence, the coefficients actually all lie in $\Z[q,q^{-1}]$, 
so property (3) holds.
Now property (3) immediately implies that bar is
continuous, so it extends uniquely to a continuous
antilinear map $-:\TT^n \rightarrow \TT^n$ still satisfying (3).
The argument in \cite[Lemma 24.1.2]{Lubook} shows that
property (2) is satisfied. Lemma~\ref{tprops}(i) gives 
that bar is an involution, whence property (1) holds, while
property (4) follows from Lemma~\ref{tprops}(iv).
Note finally that as in \cite[$\S$27.3.6]{Lubook}, the definition 
is independent of the initial choices of $n_1, n_2$.

\begin{Example}\label{bi}\rm
The bar involution on $\TT^2$ is uniquely determined 
by the following formulae.
\begin{align*}
\overline{v_a \otimes v_b} &= v_a \otimes v_b&\hbox{\!\!\!\!\!\!\!($a \leq b, a+b \neq 0$)}\\
\overline{v_a \otimes v_b} &= v_a \otimes v_b + (q^2 - q^{-2}) v_b \otimes v_a
&\hbox{\!\!\!\!\!\!\!($a > b, a+b \neq 0$)}
\end{align*}\begin{align*}
\overline{v_0 \otimes v_0} &=
v_0 \otimes v_0 + \sum_{b < 0} (q+q^{-1})(q^2 - q^{-2}) (-q^2)^{b+1} v_{b} \otimes v_{-b}\\
\overline{v_{-a} \otimes v_a} &=
v_{-a} \otimes v_a + \sum_{b > a} (q^2 - q^{-2}) (-q^2)^{a+1-b} v_{-b} \otimes v_b&\hbox{($a \geq 1$)}\\
\overline{v_{a} \otimes v_{-a}}
&=v_a \otimes v_{-a} + q^2 (q^2 - q^{-2}) v_{-a} \otimes v_a
&\hbox{($a \geq 1$)}\\
&\qquad+ \sum_{0 < b < a} (q^2 - q^{-2}) (-q^2)^{b+1-a} v_b \otimes v_{-b}\\
&\qquad+ (q - q^{-1}) (-q^2)^{1-a} v_0 \otimes v_{0}\\
&\qquad+\sum_{b < 0} q^2(q^2 - q^{-2}) (-q^2)^{b+1-a} v_{b} \otimes v_{-b}
\end{align*}
\end{Example}

Now that the bar involution has been defined on $\TT^n$, we can define
a new bilinear form $\langle.,.\rangle$ on $\TT^n$ by setting
\begin{equation}\label{lad}
\langle u,v \rangle = (u, \sigma(\overline{v}))
\end{equation}
for all $u,v \in \TT^n$,
where $(.,.)$ and $\sigma$ are as in Lemma~\ref{tricks}.
Note this makes sense even though the map $\sigma$ and the form
$(.,.)$ are not defined on the completion.

\begin{Lemma}\label{sym} $\langle.,.\rangle$ is a symmetric
bilinear form with
$\langle x u, v \rangle = \langle u, \sigma(x) v\rangle$
for all $x \in \U, u,v \in \TT^n$.
\end{Lemma}

\begin{proof}
For the second part of the lemma, 
we calculate using Lemma~\ref{tricks} to get that
\begin{align*}
\langle x u, v \rangle &= 
(x u, \sigma(\overline{v}))
=
(u, \tau(x) \sigma(\overline{v}))=
(u, \tau(\overline{\sigma\overline{(\sigma(x))}}) 
\sigma(\overline{v}))\\
&=
(u, 
\sigma(\overline{\sigma(x)}\overline{v}))
=
(u, \sigma(\overline{\sigma(x) v }))
=
\langle u, \sigma(x) v \rangle.
\end{align*}
Now let us show by induction on $n$ that $\langle.,.\rangle$ is
a symmetric bilinear form, this being obvious in case $n = 1$.
For $n > 1$, write $n = n_1 + n_2$ for $n_1,n_2 \geq 1$.
Take $u_1 \otimes u_2, v_1 \otimes v_2 \in \T^{n_1} \otimes \T^{n_2}$.
Write
$\Theta = \sum_{i \in I} x_i \otimes y_i
\in (\U \otimes \U)^\wedge$.
Recall by Lemma~\ref{tprops}(iii) that
$$
\sum_{i \in I} \sigma(\overline{x_i}) \otimes \sigma(\overline{y_i})
=
\sum_{i \in I} \overline{x_i} \otimes \overline{y_i}.
$$
Combining this with the inductive hypothesis, we calculate from
the definition of the bar involution on $\TT^n$:
\begin{align*}
\langle u_1 \otimes u_2, v_1 \otimes v_2 \rangle
&=
(u_1 \otimes u_2, \sigma(\Theta_{n_1,n_2}(\overline{v_1} \otimes \overline{v_2})))=
\sum_{i \in I}
(u_1 \otimes u_2, \sigma(x_i \overline{v_1} \otimes y_i \overline{v_2}))\\
&=
\sum_{i \in I}(u_1, \sigma(x_i \overline{v_1}))(u_2, \sigma(y_i \overline{v_2}))= \sum_{i \in I}\langle u_1, \overline{x_i} v_1 \rangle
\langle u_2, \overline{y_i} v_2 \rangle\\
&
= 
\sum_{i \in I}\langle \overline{x_i} v_1, u_1 \rangle
\langle \overline{y_i} v_2, u_2 \rangle=
\sum_{i \in I}\langle v_1, \sigma(\overline{x_i} )u_1 \rangle
\langle v_2, \sigma(\overline{y_i} )u_2 \rangle\\
&=
\sum_{i \in I}\langle v_1, \overline{x_i} u_1 \rangle
\langle v_2, \overline{y_i} u_2 \rangle= 
\langle v_1 \otimes v_2, u_1 \otimes u_2 \rangle.
\end{align*}
Hence $\langle.,.\rangle$ is symmetric.
\end{proof}

\Point{Canonical basis of $\TT^n$}
Now that we have constructed the bar involution on $\TT^n$
satisfying the properties (1)--(4) above, 
we get the following
theorem by general principles, cf. the proof of \cite[Theorem 2.17]{KL}.

\begin{Theorem}\label{thma} There exist unique topological bases
$\{T_\lambda\}_{\lambda \in \Z^n}$ and $\{L_\lambda\}_{\lambda \in \Z^n}$
for $\TT^n$ such that
\begin{itemize}
\item[(i)]
$\overline{T_\lambda} = T_\lambda$ and $\overline{L_\lambda} = L_\lambda$;
\item[(ii)]
$\overline{T_\lambda} \in N_\lambda + 
\widehat\sum_{\mu \in \Z^n}
q \Z[q] N_\mu$ and
$\overline{L_\lambda} \in M_\lambda + 
\widehat\sum_{\mu \in \Z^n}
q^{-1} \Z[q^{-1}] M_\mu$,
\end{itemize}
for all $\lambda \in \Z^n$. Actually, we have that
$\overline{T_\lambda} \in N_\lambda + 
\widehat\sum_{\mu \prec\lambda}
q \Z[q] N_\mu$ and that
$\overline{L_\lambda} \in M_\lambda + 
\widehat\sum_{\mu \prec\la}
q^{-1} \Z[q^{-1}] M_\mu$.
Also, $\omega(T_\lambda) = T_{-w_0\lambda}$
and $\omega(L_\lambda) = L_{-w_0\lambda}$.
\end{Theorem}

\begin{Example}\label{t2}\rm
Suppose that $n = 2$.
Using Example~\ref{bi}, one checks:
\begin{align*}
T_{(a,b)} &= N_{(a,b)} &\!\!\!\!\!(a  \leq b, a+b\neq 0)\\
T_{(a,b)} &= N_{(a,b)}+q^2 N_{(b,a)} &\!\!\!\!\!(a  > b, a+b\neq 0)\\
T_{(-a,a)} &= N_{(-a,a)} + q^2 N_{(-a-1,a+1)}&(a \geq 1)\\
T_{(a,-a)} &= N_{(a,-a)} + q^2( N_{(a-1,1-a)}+N_{(1-a,a-1)})
+q^4 N_{(-a,a)}&(a \geq 2)\\
T_{(0,0)} &= N_{(0,0)} + (q+q^3) N_{(-1,1)}\\
T_{(1,-1)} &= N_{(1,-1)} + q N_{(0,0)} + q^4 N_{(-1,1)}
\end{align*}
Note in this example that each $T_\lambda$ is a finite sum of $N_\mu$'s.
I conjecture that this is true in general.
On the other hand, the $L_\lambda$'s need not be finite sums
of $M_\mu$'s even for $n = 2$.
\end{Example}

We call the topological basis $\{T_\lambda\}_{\lambda \in \Z^n}$ the
{\em canonical basis} of $\TT^n$ and $\{L_\lambda\}_{\lambda \in \Z^n}$
the {\em dual canonical basis}. 
Let us introduce notation for the coefficients: let
\begin{equation}\label{tde}
T_\lambda = \sum_{\mu \in \Z^{n}} t_{\mu,\lambda}(q) N_\mu,
\qquad
L_\lambda = \sum_{\mu \in \Z^{n}} l_{\mu,\lambda}(q) M_\mu
\end{equation}
for polynomials $t_{\mu,\lambda}(q) \in \Z[q]$ and 
$l_{\mu,\lambda}(q) \in \Z[q^{-1}]$.
We know that $t_{\mu,\lambda}(q) = l_{\mu,\lambda}(q) = 0$ unless
$\mu \preceq \lambda$, and that
$t_{\lambda,\lambda}(q) = l_{\la,\la}(q) = 1$.

\begin{Lemma}\label{thmd}
For $\lambda, \mu \in \Z^n$,
$\langle L_\lambda, T_{-\mu} \rangle = 
\delta_{\lambda,\mu}$.
\end{Lemma}

\begin{proof}
A calculation using the definition of $\langle.,.\rangle$ shows that
\begin{align}\label{fa}
\langle L_\lambda, T_{-\mu} \rangle&=
\sum_{\mu \preceq \nu \preceq\lambda}
l_{\nu,\lambda}(q) t_{-\nu,-\mu}(q^{-1}),\\
\langle T_{-\mu},L_\lambda \rangle
&=
\sum_{\mu \preceq \nu \preceq\lambda}
l_{\nu,\lambda}(q^{-1}) t_{-\nu,-\mu}(q).
\end{align}
Hence, $\langle L_\lambda,T_{-\mu}\rangle$
equals $1$ if $\lambda = \mu$ and belongs to
$q^{-1} \Z[q^{-1}]$ if $\lambda \neq \mu$.
Similarly,
$\langle T_{-\mu},L_\lambda\rangle$  equals $1$ if 
$\lambda = \mu$ and belongs to
$q \Z[q]$ if $\lambda \neq \mu$.
But $\langle L_\lambda,T_{-\mu}\rangle=\langle T_{-\mu},L_\lambda\rangle$
by Lemma~\ref{sym}.
\end{proof}

\begin{Corollary}\label{fne}
For $\lambda \in \Z^n$,
$M_\lambda = \sum_{\mu \in\Z^n} t_{-\lambda,-\mu}(q^{-1}) L_\mu$.
\end{Corollary}

\begin{proof}
By Lemma~\ref{thmd}, we can write
$M_\la = \sum_{\mu \in \Z^n} \langle M_\la,T_{-\mu} \rangle L_\mu$.
Now a calculation from the definition of the form $\langle.,.\rangle$
gives that
$\langle M_\la, T_{-\mu} \rangle = t_{-\la,-\mu}(q^{-1})$.
\end{proof}

\begin{Example}\label{m2}\rm
Suppose that $n = 2$.
Using Example~\ref{t2}, one checks:
\begin{align*}
M_{(a,b)} &= L_{(a,b)} &\!\!\!\!\!\!\!\!\!\!\!\!(a  \leq b, a+b\neq 0)\\
M_{(a,b)} &= L_{(a,b)} + q^{-2} L_{(b,a)}&\!\!\!\!\!\!\!\!\!\!\!\!(a  > b, a+b\neq 0)\\
M_{(-a,a)} &= L_{(-a,a)} + q^{-2} L_{(-a-1,a+1)}&(a \geq 1)\\
M_{(a,-a)} &= L_{(a,-a)} + q^{-2}( L_{(a-1,1-a)}+L_{(-a-1,a+1)})
+q^{-4} L_{(-a,a)}&(a \geq 2)\\
M_{(0,0)} &= L_{(0,0)} + q^{-1} L_{(-1,1)}\\
M_{(1,-1)} &= L_{(1,-1)} + (q^{-1}+q^{-3}) L_{(0,0)} +   q^{-2} L_{(-2,2)}+q^{-4} L_{(-1,1)}\!\!\!\!
\end{align*}
I conjecture for arbitrary $n$ 
that each $M_\la$ is always a finite linear combination
of $L_\mu$'s.
\end{Example}

\Point{Crystal structure}\label{crys}
Now we describe the crystal structure underlying the module 
$\T^n$. The basic reference followed here is \cite{Ka}.
Let $\A$ be the subring of $\Q(q)$ consisting of rational functions
having no pole at $q = 0$.
Evaluation at $q = 0$ induces an isomorphism $\A / q \A \rightarrow \Q$.

Let $\V_\A$ be the $\A$-lattice in $\V$ spanned by the $v_a$'s.
Then, $\V_\A$ together with the basis of the $\Q$-vector space
$\V_\A / q \V_\A$
given by the images of the $v_a$'s is a lower {crystal base}
for $\V$ at $q = 0$ in the sense of \cite[4.1]{Ka}.
Write $\tilde E_i, \tilde F_i$ for the corresponding crystal 
operators. Rather than view these as operators on the
crystal base $\{v_a + q \V_\A\}_{a \in \Z}$, we will view them simply
as operators on the set $\Z$ parameterizing the
crystal base. Then, the crystal graph is as follows:
$$\dots\longrightarrow
-3 \stackrel{\tilde F_2}{\longrightarrow}
-2 \stackrel{\tilde F_1}{\longrightarrow}
-1 \stackrel{\tilde F_0}{\longrightarrow}
0 \stackrel{\tilde F_0}{\longrightarrow}
1 \stackrel{\tilde F_1}{\longrightarrow}
2 \stackrel{\tilde F_2}{\longrightarrow}
3 \longrightarrow\dots.
$$
Thus, $\tilde F_i (a)$ equals $a+1$ if $a = i$ or $a = -i-1$,
$\varnothing$ otherwise, and 
$\tilde E_i (a)$ equals $a-1$ if $a = i+1$ or $a = -i$,
$\varnothing$ otherwise.
The maps $\eps_i, \phi_i:\Z \rightarrow \N$ defined by
$$
\eps_i(a) = \max\{k \geq 0\:|\:\tilde E_i^k a \neq \varnothing\},
\qquad
\phi_i(a) = \max\{k \geq 0\:|\:\tilde F_i^k a \neq \varnothing\}
$$
only ever take the values $0, 1$ or $2$ (the last possibility occurring only 
if $i = 0$).

Since $\T^n$ is a tensor product of $n$ copies of $\V$, it has
an induced crystal structure.
The crystal lattice is $\T^n_\A$, namely, the $\A$-span of the basis
$\{N_\lambda\}_{\lambda \in \Z^n}$, and the images of the
$N_\lambda$'s in $\T_\A^n / q\T^n_\A$ give the crystal base.
Like in the previous paragraph, we will view the crystal
operators as maps on the underlying set $\Z^n$ parametrizing the
crystal base. However, 
we will denote them by $\tilde E'_i, \tilde F'_i$, 
since we want to reserve the unprimed symbols for something else later on.
In order to describe them explicitly, we introduce a little more
combinatorial notation.
Given $\lambda \in \Z^n$ and $i \geq 0$, 
let $(\sigma_1,\dots,\sigma_n)$ be the {\em $i$-signature}
of $\lambda$, namely, the sequence defined by
\begin{equation}\label{isig}
\sigma_r = \left\{
\begin{array}{ll}
+&\hbox{if $i \neq 0$ and $\lambda_r = i$ or $-i-1$,}\\
-&\hbox{if $i \neq 0$ and $\lambda_i = i+1$ or $-i$,}\\
++&\hbox{if $i = 0$ and $\lambda_r = -1$,}\\
-+&\hbox{if $i = 0$ and $\lambda_r = 0$,}\\
--&\hbox{if $i = 0$ and $\lambda_r = 1$,}\\
0&\hbox{otherwise.}
\end{array}\right.
\end{equation}
Form the {\em reduced $i$-signature} by successively replacing
subwords of $(\sigma_1,\dots,\sigma_n)$
of the form $+-$ (possibly separated by $0$'s in between)
with $0$'s until we are left with 
a sequence of the form
$(\tilde \sigma_1, \dots, \tilde \sigma_n)$ 
in which no $-$ appears after a $+$.
For $r = 1,\dots,n$, let $\delta_r$ be the $n$-tuple
with a $1$ in the $r$th position and $0$'s elsewhere.
Then,
\begin{align*}
\tilde E'_i (\lambda) &= \left\{
\begin{array}{ll}
\varnothing&\hbox{if there are no $-$'s in the reduced $i$-signature,}\\
\lambda - \delta_r&\hbox{otherwise, where the rightmost 
$-$ occurs in $\tilde \sigma_r$,}
\end{array}
\right.\\
\tilde F'_i (\lambda) &= \left\{
\begin{array}{ll}
\varnothing&\hbox{if there are no $+$'s in the reduced $i$-signature,}\\
\lambda+\delta_r&\hbox{otherwise, where the leftmost $+$ occurs in $\tilde \sigma_r$,}\\
\end{array}
\right.\\
\eps'_i(\lambda) &= 
\hbox{the total number of $-$'s in the reduced $i$-signature},\\
\phi'_i(\lambda) &=
\hbox{the total number of $+$'s in the reduced $i$-signature}.
\end{align*}
We have now described the crystal 
$(\Z^{n}, \tilde E'_i, \tilde F'_i, \eps'_i, \phi'_i, \wt)$
associated to  the module $\T^{n}$ purely combinatorially.

\begin{Example}\rm\label{eg1}
Consider
$\lambda = (1,2,0,-3,-2,-1,0,1) \in \Z^8$.
The $1$-signature is $(+,-,0,0,+,-,0,+)$.
Hence the reduced $1$-signature is $(0,0,0,0,0,0,0,+)$, so 
we get that $\tilde E'_1 \lambda = \varnothing$ and
$\tilde F_1' \lambda = \lambda + \delta_8$.
On the other hand, the $0$-signature is
$(--,0,-+,0,0,++,-+,--)$, which reduces
to
$(--,0,-+,0,0,0,0,0)$, so 
$\tilde E'_0 \lambda = \lambda - \delta_3$,
$\tilde F'_0 \lambda = \lambda + \delta_3$.
\end{Example}

The following lemma is a general property of canonical bases/lower 
global crystal bases. It follows 
ultimately from \cite[Proposition 5.3.1]{KaG}.
See \cite[Theorem 2.31]{KL} for a similar situation.

\begin{Lemma}\label{thmc}
Let $\lambda \in \Z^{n}$ and $i \geq 0$.
\begin{itemize}
\item[(i)]
$E_i T_\lambda = 
[\phi'_i(\lambda)+1]_i T_{\tilde E'_i(\lambda)} + 
\sum_{\mu \in \Z^{n}}
u_{\mu,\lambda}^i T_\mu$
where $u_{\mu,\lambda}^i\in
qq_i^{1 - \phi'_i(\mu)} \Z[q]$
is zero unless
$\eps'_j(\mu) \geq \eps'_j(\lambda)$ for all $j \geq 0$.
\item[(ii)]
$F_i T_\lambda = [\eps'_i(\lambda)+1]_i T_{\tilde F'_i(\lambda)} + 
\sum_{\mu \in \Z^{n}}
v_{\mu,\lambda}^i T_\mu$
where $v_{\mu,\lambda}^i\in qq_i^{1 - \eps'_i(\mu)} \Z[q]$
is zero unless
$\phi'_j(\mu) \geq \phi'_j(\lambda)$ for all $j \geq 0$.
\end{itemize}
(In (i) resp. (ii), the first term on the right hand side 
should be omitted if
$\tilde E'_i(\la)$ resp. 
$\tilde F'_i(\la)$ is $\varnothing$.)
\end{Lemma}

Motivated by Lemmas~\ref{sym} and \ref{thmd}, we also
introduce the {\em dual crystal operators} defined by
\begin{align}
\tilde E_i^*(\lambda) = -\tilde F'_{i} (-\lambda),\qquad
&\tilde F_i^*(\lambda) = -\tilde E'_{i} (-\lambda),\\
\eps_i^*(\lambda) =  \phi'_{i}(-\lambda),
\qquad
&\phi_i^*(\lambda) =  \eps'_{i}(-\lambda).
\end{align}
These can be described explicitly in a similar way to the above:
for fixed  $i \geq 0$ and $\la \in \Z^{n}$, let
$(\sigma_1,\dots,\sigma_n)$ be the $i$-signature of $\la$ defined according
to (\ref{isig}). First, replace all $\sigma_r$ that equal $-+$
with $+-$.
Now form the dual reduced $i$-signature $(\tilde \sigma_1,\dots,
\tilde \sigma_n)$ from this by repeatedly replacing
subwords of $(\sigma_1,\dots,\sigma_n)$ 
of the form  $-+$ (possibly separated by $0$'s)
by $0$'s, until no $+$ appears after a $-$.
Finally, 
\begin{align*}
\tilde E^*_i (\lambda) &= \left\{
\begin{array}{ll}
0&\hbox{if there are no $-$'s in the dual reduced $i$-signature,}\\
\lambda-\delta_r&\hbox{otherwise, where the leftmost $-$ occurs in
$\tilde \sigma_r$,}
\end{array}
\right.\\
\tilde F_i^* (f) &= \left\{
\begin{array}{ll}
0&\hbox{if there are no $+$'s in the dual reduced $a$-signature,}\\
\la+\delta_r&\hbox{otherwise, where the rightmost $+$ occurs in $\tilde\sigma_r$,}
\end{array}
\right.\\
\eps^*_i(\lambda) &=
\hbox{the total number of $-$'s in the dual reduced $i$-signature},\\
\phi_i^*(\lambda) &=
\hbox{the total number of $+$'s in the dual reduced $i$-signature}.
\end{align*}
In this way, we obtain the dual crystal structure
$(\Z^{n}, \tilde E^*_i, \tilde F^*_i, \eps_i^*, 
\phi_i^*, \wt)$.

\begin{Lemma}\label{cthm}
Let $\lambda \in \Z^{n}$ and $i \geq 0$.
\begin{itemize}
\item[(i)]
$E_i L_\lambda = [\eps_i^*(\la)]_i L_{\tilde E_i^* (\la)} + 
\sum_{\mu \in \Z^{n}}
w_{\mu,\la}^i L_\mu$
where $w_{\mu,\la}^i\in qq_i^{1 - \eps_i^*(\la)} \Z[q]$
is zero unless
$\phi_j^*(\mu) \leq \phi_j^*(\la)$ for all $j \geq 0$.
\item[(ii)]
$F_i L_\la = [\phi_i^*(\la)]_i L_{\tilde F_i^* (\la)} + 
\sum_{\mu \in \Z^{n}}
x_{\mu,\la}^i L_\mu$
where $x_{\mu,\la}^i\in qq_i^{1 - \phi_i^*(\la)} \Z[q]$
is zero unless
$\eps_j^*(\mu) \leq \eps_j^*(\la)$ for all $j \geq 0$.
\end{itemize}
\end{Lemma}

\begin{proof}
Dualize Lemma~\ref{thmc} using Lemmas~\ref{sym} and \ref{thmd}. 
\end{proof}

\begin{Remark}\label{fness}\rm
(i) Suppose we are given $\eps_i,\phi_i\in\N$ for all $i \geq 0$.
One can show from the combinatorial description of
the maps $\eps_i^*, \phi_i^*$ above that
there exist only finitely many
$\la \in \Z^{n}$ with $\eps^*_i(\la) = \eps_i$ and 
$\phi^*_i(\la) = \phi_i$
for all $i \geq 0$.

(ii) Using (i), one deduces easily that 
all but finitely many terms of the sums occurring in Lemma~\ref{cthm} 
are zero, i.e. both $E_i L_\lambda$ and $F_i L_\lambda$ are finite
$\Z[q,q^{-1}]$-linear combinations
of $L_\mu$'s.  We will not make us of this observation in the remainder of the
article.

(iii) If the finiteness conjecture made in Example~\ref{t2} holds,
then it is also the case that
all but finitely many terms of the sums occurring in Lemma~\ref{thmc} are
zero, i.e.
both $E_i T_\lambda$ and $F_i T_\lambda$ are finite
$\Z[q,q^{-1}]$-linear combinations
of $T_\mu$'s, cf. the argument in 
the last paragraph of the proof of
Lemma~\ref{argb} below.
\end{Remark}

\section{Exterior algebra}\label{ea}

\Point{Exterior powers}\label{ep}
Define $\K$ to be the two-sided ideal of the tensor algebra $\T$
generated by the
vectors
$$
\begin{array}{lr}
v_a \otimes v_a&(a \neq 0)\\
v_a \otimes v_b + q^{2} v_b \otimes v_a&(a > b, a +b \neq 0)\\
v_{a}\otimes v_{-a}
+q^2(v_{a-1} \otimes v_{1-a}+v_{1-a} \otimes v_{a-1})
+ q^4 v_{-a} \otimes v_{a}
&(a \geq 2)\\
v_{1} \otimes v_{-1}+ q v_{0} \otimes v_0+ q^4 v_{-1} \otimes v_{1},
\end{array}
$$
for all admissible $a,b \in \Z$.
These relations are a limiting case of 
the relations in \cite[Proposition 2.3]{JMO}
(with $q$ replaced by $q^2$), hence $\K$ is invariant under the action of $\U$.
Let $\F := \T / \K$.
Since $\K = \bigoplus_{n \geq 0} \K^n$ 
is a homogeneous ideal of $\T$, 
$\F$ is also 
graded as $\F = \bigoplus_{n \geq 0} \F^n$, where $\F^n = \T^n / \K^n$.
We view the space $\F^n$ as a quantum analogue
of the exterior power $\bigwedge^n \V$ in type $\mathfrak{b}_{\infty}$.
As usual, we will write $u_1 \wedge \dots \wedge u_n$ for the
image of $u_1 \otimes \dots \otimes u_n  \in \T^n$ under the quotient map
$\pi:\T^n \rightarrow \F^n$.

As in the introduction, let $\Z^n_+$ denote the set of all tuples
$\lambda \in \Z^n$ such that
$\lambda_r > \lambda_{r+1}$ if $\lambda_r \neq 0$,
$\lambda_r \geq \lambda_{r+1}$ if $\lambda_r = 0$,
for each $r = 1,\dots,n-1$.
For $\lambda \in \Z^n_+$, let
\begin{equation}
F_\lambda = \pi(N_{w_0 \lambda}) = 
v_{\lambda_n} \wedge \dots \wedge v_{\lambda_1} \in \F^n.
\end{equation}
The following lemma follows  from the defining relations for $\K^n$.

\begin{Lemma}\label{handy} For $\lambda \in \Z^n$, $\pi(N_{w_0\lambda})$ 
equals $F_{\lambda}$ if $\lambda \in \Z_n^+$,
otherwise $\pi(N_{w_0\lambda})$ is a $q \Z[q]$-linear combination of
$F_\mu$'s for $\mu \in \Z^n_+$ with $\mu \succeq \lambda$.
\end{Lemma}

This shows that the elements 
$\{F_\lambda\}_{\lambda \in \Z^n_+}$ span
$\F^n$. In fact, one can check routinely using Bergman's diamond lemma
\cite[1.2]{Berg} that:

\begin{Lemma}\label{bt}
The vectors $\{F_\lambda\}_{\lambda \in \Z^n_+}$ give a basis for
$\F^n$.
\end{Lemma}

Let $\KK^n$ be the closure of $\K^n$ in $\TT^n$ and 
$\FF^n := \TT^n / \KK^n$, giving 
a completion of the vector space $\F^n$.
The vectors $\{F_\lambda\}_{\lambda \in \Z^n_+}$ give a topological basis
for $\FF^n$.
Note that $\sigma:\T^n\rightarrow\T^n$
leaves $\K^n$ invariant, hence induces
$\sigma:\F^n\rightarrow\F^n$.
Similarly the continuous automorphism $\omega:\TT^n\rightarrow\TT^n$
leaves $\KK^n$ invariant, so induces
$\omega:\FF^n\rightarrow\FF^n$ with 
$\omega(F_\lambda) = F_{-w_0\lambda}$.

\Point{\boldmath Canonical basis of $\FF^n$}
Now we construct the canonical basis of $\FF^n$.
To start with, we need a bar involution.

\begin{Lemma}\label{bi2} The bar involution on $\TT^n$ leaves $\KK^n$ invariant,
hence induces a continuous antilinear involution
$-:\FF^n \rightarrow \FF^n$ such that
\begin{itemize}
\item[(1)] $\overline{xv} = \overline{x}\,\overline{v}$ for all
$x \in \U, v \in \FF^n$;
\item[(2)] $\overline{F_\lambda} \in F_\lambda + 
\sum_{\mu\succ\lambda} \Z[q,q^{-1}]F_\mu$ for all $\lambda \in \Z^n_+$;
\item[(3)] $\overline{\omega(v)} = \omega(\overline{v})$ for all
$v\in\FF^n$.
\end{itemize}
\end{Lemma}

\begin{proof}
In the case $n = 2$, all the generators of $\K^2$ are bar invariant
by Example~\ref{t2}, hence $\K^2$ is bar invariant.
In general, $\K^n$ is spanned by vectors of the form
$v \otimes k \otimes w$ for $v \in \T^{n_1}, k \in \K^2,
w \in T^{n_2}$ and some $n_1,n_2 \geq 0$ with $n_1+n_2+2 = n$.
By the definition of the bar involution,
$$
\overline{v \otimes k \otimes w}
=
\Theta_{n_1+2,n_2} (\overline{v} \otimes \Theta_{2,n_2} (\overline{k} \otimes \overline{w})).
$$
We have already shown that $\overline{k} \in \K^2$, and $\K^2$ is $\U$-invariant, hence this belongs to $\KK^n$.
This shows that $\overline{\K^n} \subset \KK^n$, hence $\KK^n$ itself is 
bar invariant by continuity.
Now properties (1) and (3) are immediate from the analogous properties 
of the bar involution on $\TT^n$.
For property (2), take $\lambda \in \Z^n_+$. We know
$\overline{N_{w_0 \lambda}}$ equals $N_{w_0\lambda}$
plus {a $\Z[q,q^{-1}]$-linear 
combination of $N_{w_0 \mu}$'s
with $w_0  \mu \prec w_0 \lambda$}, or 
equivalently, $\mu \succ \lambda$.
By Lemma~\ref{handy}, $\pi(N_{w_0 \mu})$ is a
$\Z[q,q^{-1}]$-linear combination of $F_\nu$'s
with $\nu \succeq \mu$.
Hence (2) holds.
\end{proof}

We get the following theorem by general principles, just
as in Theorem~\ref{thma} before.

\begin{Theorem}\label{thmb} There exists a unique topological basis
$\{U_\lambda\}_{\lambda \in \Z^n_+}$ for $\FF^n$ such that
$\overline{U_\lambda} = U_\lambda$
and ${U_\lambda} \in F_\lambda + 
\widehat\sum_{\mu \in \Z_+^n}
q \Z[q] F_\mu$,
for all $\lambda \in \Z_+^n$. Actually, we have that
${U_\lambda} \in F_\lambda + 
\widehat\sum_{\mu \succ\lambda}
q \Z[q] F_\mu.$
Also, $\omega(U_\la) = U_{-w_0\la}$.
\end{Theorem}

\begin{Example}\rm\label{ex1}
For $n = 2$, one deduces from Example~\ref{t2} and Lemma~\ref{tss} below that:
\begin{align*}
U_{(a,b)} &= F_{(a,b)} &(a > b, a+b \neq 0)\\
U_{(a,-a)} &= F_{(a,-a)} + q^2 F_{(a+1,-a-1)}&(a \geq 1)\\
U_{(0,0)} &= F_{(0,0)} + (q+q^3) F_{(1,-1)}
\end{align*}
\end{Example}

We call the topological basis 
$\{U_\lambda\}_{\lambda \in \Z^n_+}$ the {\em canonical basis}
of $\FF^n$. Write
\begin{equation}\label{ude}
U_\lambda = \sum_{\mu \in \Z^{n}_+} u_{\mu,\lambda}(q) F_\mu
\end{equation}
for polynomials $u_{\mu,\lambda}(q) \in \Z[q]$.
We know that $u_{\mu,\lambda}(q) = 0$ unless
$\mu \succeq \lambda$, and that
$u_{\lambda,\lambda}(q) =1$. The following lemma explains the relationship
between $U_\lambda$ and the canonical basis element
$T_\lambda$ of $\TT^n$ constructed earlier.

\begin{Lemma}\label{tss} For $\lambda \in \Z^n$, we have that
$$
\pi(T_{w_0 \lambda}) = \left\{\begin{array}{ll}
U_{\lambda}&\hbox{if $\lambda \in \Z^n_+$,}\\
0&\hbox{if $\lambda \notin \Z^n_+$.}
\end{array}\right.
$$
\end{Lemma}

\begin{proof}
Suppose first that $\lambda \in \Z^n_+$.
We know that $T_{w_0\lambda}$ equals 
$N_{w_0 \lambda}$ plus
a $q\Z[q]$-linear combination of $N_{w_0 \mu}$'s
with $\mu \succ \lambda$.
Moreover, 
by Lemma~\ref{handy}, $\pi(N_{w_0 \mu})$ is a
$\Z[q]$-linear combination of $F_{\nu}$'s with
$\nu \succeq \mu$. Hence,
$$
\pi(T_{w_0 \lambda})
= F_{\lambda} + 
(\hbox{a $q\Z[q]$-linear combination of $F_{\mu}$'s
with $\mu \succ \lambda$}).
$$
Since it is automatically bar invariant, it must equal $U_\lambda$
by the uniqueness in Theorem~\ref{thmb}.
An entirely similar argument in case $\lambda \notin \Z^n_+$
shows that
\begin{align*}
\pi(T_{w_0 \lambda}) &= 
(\hbox{a $q \Z[q]$-linear combination of $F_{\mu}$'s
with $\mu \succ \lambda$})\\
&= 
(\hbox{a $q \Z[q]$-linear combination of $U_{\mu}$'s
with $\mu \succ \lambda$}).
\end{align*}
Since it is bar invariant, it must be zero.
\end{proof}

\begin{Corollary}\label{tbase}
The vectors $\{T_{w_0\lambda}\}_{\lambda \notin \Z^n_+}$
form a topological basis for $\KK^n$.
\end{Corollary}

\Point{Dual canonical basis}\label{dcb}
For $\lambda \in \Z^n_+$, define
\begin{equation}\label{edef}
E_\lambda =
\sum_{\mu \in\Z^n_+}
u_{-w_0\lambda,-w_0\mu}(q^{-1}) L_\mu.
\end{equation}
Let $\E^n$ be the subspace of $\TT^n$ 
spanned by the $\{E_\lambda\}_{\lambda \in \Z^n_+}$.
Since there are only finitely many $\mu \in \Z^n_+$ with $\mu \preceq \lambda$,
we see that $E_\lambda$ is a finite linear combination of $L_\mu$'s,
and vice versa.
So the vectors $\{L_\lambda\}_{\lambda \in \Z^n_+}$ also form a basis
for $\E^n$.

\begin{Example}\rm For $n = 2$, we have by Example~\ref{ex1}
that:
\begin{align*}
E_{(a,b)} &= L_{(a,b)} &(a > b, a+b \neq 0)\\
E_{(a,-a)} &= L_{(a,-a)} + q^{-2} L_{(a-1,-a+1)}&(a \geq 2)\\
E_{(1,-1)} &= L_{(1,-1)} + (q^{-1}+q^{-3}) L_{(0,0)}\\
E_{(0,0)} &= L_{(0,0)}
\end{align*}
\end{Example}

Let $\EE^n$ be the closure of $\E^n$ in $\TT^n$.
By Corollary~\ref{tbase} and Lemma~\ref{thmd},
$\EE^n$ and $\KK^n$ are orthogonal with respect to the bilinear
form $\langle.,.\rangle$. 
Hence we get induced a well-defined pairing
$\langle.,.\rangle$ between $\EE^n$ and $\FF^n$.
By Lemmas~\ref{tss} and \ref{thmd}, we have that
\begin{equation}\label{dty}
\langle L_{-w_0\lambda}, U_{\mu}\rangle = 
\delta_{\lambda,\mu}
\end{equation}
for all $\lambda,\mu \in \Z^n_+$.

\begin{Lemma}\label{sub} $\EE^n$ is a $\U$-submodule
of $\TT^n$.
\end{Lemma}

\begin{proof}
It suffices to show that $E_i L_\lambda$ and $F_i L_\lambda$
both belong to $\EE^n$ for each $\lambda \in \Z^n_+$ and $i \geq 0$.
Write
$E_i L_\lambda = \sum_{\mu \in \Z^m} b_{\mu,\lambda}(q) L_\mu.$
Apply $\langle.,T_{-\mu}\rangle$ to both sides and use
Lemma~\ref{thmd} to get
$b_{\mu,\lambda}(q) = \langle L_\lambda, E_i T_{-\mu}\rangle.$
For $\mu \notin \Z_+^n$, $T_{-\mu}$ belongs to $\KK^n$.
But $\KK^n$ is $\U$-invariant, hence $E_i T_{-\mu}$ belongs to $\KK^n$.
So we get that $b_{\mu,\lambda}(q) = 0$ for all $\mu\notin\Z^n_+$
by Corollary~\ref{tbase}.
Hence $E_i L_\lambda$ belongs to $\EE^n$ still, and similarly
for $F_i L_\lambda$.
\end{proof}

The bilinear form $(.,.)$ on $\T^n$ also induces a pairing
$(.,.)$ between $\E^n$ and $\F^n$, actually by (\ref{lad}) 
we have that
\begin{equation}\label{lad2}
(u,v) = \langle u, \overline{\sigma(v)}\rangle
\end{equation} 
for each $u\in\E^n,v\in\F^n$.

\begin{Lemma}\label{bf}
For all $\lambda, \mu \in \Z^n_+$, we have that
$(E_\lambda, N_{\mu}) = q^{z(\lambda)}\delta_{\lambda,\mu}.$
\end{Lemma}

\begin{proof}
Let $F_{\mu} = \sum_{\gamma \in \Z^n_+} v_{\gamma,\mu}(q) U_{\gamma}$, so
$\sum_{\gamma \in \Z^n_+} u_{\mu,\gamma}(q) v_{\gamma,\nu}(q) = 
\delta_{\mu,\nu}.$
Now calculate:
\begin{align*}
(E_{-w_0\lambda},N_{-w_0\mu}) &= 
\langle E_{-w_0\lambda}, \overline{\sigma(N_{-w_0\mu})}\rangle\\&=
q^{z(\mu)} \langle E_{-w_0\lambda}, \overline{N_{w_0\mu}}\rangle=
q^{z(\mu)} \langle E_{-w_0\lambda}, \overline{F_{\mu}} \rangle\\&=
q^{z(\mu)}
\sum_{\nu,\gamma \in \Z^n_+} 
u_{\la,\nu}(q^{-1})
v_{\gamma,\mu}(q^{-1}) \langle L_{-w_0\nu}, U_{\gamma}\rangle = q^{z(\mu)} \delta_{\la,\mu}.
\end{align*}
\end{proof}

Now we get the following characterization of the basis
$\{L_\lambda\}_{\lambda \in \Z^n_+}$ 
in terms of the restriction of the bar involution to $\EE^n$.

\begin{Theorem}
For $\lambda \in \Z^n_+$, $L_\lambda$ is the unique
element of $\EE^n$ such that $\overline{L_\lambda} = L_\lambda$
and $L_\lambda \in E_\lambda + \widehat\sum_{\mu \in \Z^n_+} q^{-1}\Z[q^{-1}] E_\mu$.
Moreover,
\begin{equation}\label{lfin}
L_\lambda = \sum_{\mu \in \Z^n_+} l_{\mu,\lambda}(q) E_\mu,
\end{equation}
where $l_{\mu,\lambda}(q)$ is as in (\ref{tde}).
\end{Theorem}

\begin{proof}
For each $\mu \in \Z^n_+$,
$(L_\lambda,N_\mu) = q^{z(\mu)} l_{\mu,\lambda}(q)$.
Now apply Lemma~\ref{bf} to deduce that
$L_\lambda = \sum_{\mu \in \Z^n_+} l_{\mu,\lambda}(q) E_\mu$.
Finally, if $L_\lambda'\in 
E_\lambda + \sum_{\mu \in \Z^n_+} q^{-1}\Z[q^{-1}] E_\mu$ 
is another bar invariant element of
$\EE^n$,
then $L_\lambda - L_\lambda'$ is bar invariant
and can be expressed as a $q^{-1}\Z[q^{-1}]$-linear combination of
$L_\nu$'s. Hence it must be zero.
\end{proof}

To state the next lemma,
we define polynomials $a_{\la,\mu}(q) \in \Z[q]$
for each $\lambda\in\Z_+^n, \mu \in \Z^n$ by
\begin{equation}\label{adef}
\pi(N_{w_0\mu}) = \sum_{\la \in \Z^n_+}
a_{\la,\mu}(q) F_{\la}.
\end{equation}
Recalling Lemma~\ref{handy}, we have that
$a_{\la,\mu}(q) = 0$ unless $\la \succeq \mu$.

\begin{Lemma}\label{ems} For each $\lambda \in \Z^n_+$,
$E_\lambda = \sum_{\mu\in\Z^n} a_{-w_0\la,-w_0\mu}(q^{-1}) M_\mu.$
\end{Lemma}

\begin{proof}
Applying the antilinear map 
$\sigma$ to (\ref{adef}) gives that
$$
\pi(N_{\mu}) = \sum_{\la \in\Z_+^n} q^{z(\mu)-z(\lambda)} a_{-w_0\la,-w_0\mu}(q^{-1})
\pi(N_{\lambda})
$$
for each $\mu \in \Z^n$.
Hence, invoking Lemma~\ref{bf},
$(E_\lambda,N_\mu) = q^{z(\mu)} a_{-w_0\la,-w_0\mu}(q^{-1})$
for each $\mu \in \Z^n$.
The lemma follows.
\end{proof}

Finally in this subsection, we describe the action
of $\U$ on the basis $\{E_\lambda\}_{\lambda\in\Z^n_+}$ of $\E^n$
explicitly. In particular, this shows that $\E^n$ itself
is a $\U$-submodule of $\TT^n$, as could also be proved using  Lemma~\ref{sub}
and Remark~\ref{fness}(ii).

\begin{Lemma}\label{expact} Let $\lambda \in \Z^n_+$ and $i \geq 0$,
and let $(\sigma_1,\dots,\sigma_n)$ be the $i$-signature
of $\lambda$ defined according to (\ref{isig}).
Then,
\begin{align*}
E_i E_\lambda &= 
q^{-(\alpha_i, \eps_{\lambda_{r+1}} + \dots + \eps_{\lambda_n})}
\!\!\!\!\sum_{\substack{r\,\text{with}\,\lambda-\delta_r\in\Z^n_+,\\\sigma_r = -,-+\,\text{or}\,--}}
c_{\la,r}(q) E_{\lambda-\delta_r},\\
F_i E_\lambda &=
q^{(\alpha_i, \eps_{\lambda_1} + \dots + \eps_{\lambda_{r-1}})}
\sum_{\substack{r\,\text{with}\,\lambda+\delta_r\in\Z^n_+,\\\sigma_r = +,-+\,\text{or}\,++}}
c_{\la,r}(q) E_{\lambda+\delta_r},
\end{align*}
where
$c_{\la,r}(q) = \left\{
\begin{array}{ll}
(q+q^{-1}) \sum_{s=0}^{z(\lambda)} (-q^{-2})^s
&\hbox{if $\sigma_r = --$ or $++$,}\\
1&\hbox{otherwise.}
\end{array}\right.$
\end{Lemma}

\begin{proof}
We  sketch  the proof for $E_i$.
By Lemmas~\ref{sub}, \ref{bf} and \ref{tricks}(ii), we may write
$E_i E_\lambda = \sum_{\mu \in \Z^n_+} c_{\mu,\lambda}(q) E_\mu$
where
$$
c_{\mu,\lambda}(q) = q^{-z(\mu)} (E_i E_\lambda, N_\mu)
=
q^{-z(\mu)} (E_\lambda, q_i F_i K_i^{-1} N_\mu).
$$
The right hand side is computed using
Lemma~\ref{bf} and the fact that $E_\lambda$ is orthogonal to $\K^n$.
\end{proof}

\Point{\boldmath Crystal structure}\label{cs}
The crystal structure underlying the canonical basis
$\{U_\lambda\}_{\lambda \in \Z^n_+}$ of $\FF^n$ is easily 
deduced from the results of \ref{crys} and Lemma~\ref{tss}.
Let us denote the resulting crystal operators on the index
set $\Z^n_+$ parameterizing the bases of $\FF^n$ by
$\tilde E_i, \tilde F_i, \eps_i, \phi_i$.
By definition,
\begin{align}
\tilde E_i(\lambda) = w_0 \tilde E'_{i} (w_0  \lambda),\qquad
&\tilde F_i(\lambda) = w_0 \tilde F'_{i} (w_0 \lambda),\\
\eps_i(\lambda) =  \eps_{i}'(w_0\lambda),
\qquad
&\phi_i(\lambda) =  \phi_{i}'(w_0 \lambda),
\end{align}
where $\tilde E_i', \tilde F_i', \eps_i'$ and $\phi_i'$ are as in
\ref{crys}.
By properties of the automorphism $\omega$ (or by directly
checking all of the cases listed below),
the operators $\tilde E_i, \tilde F_i, \eps_i, \phi_i$ 
are {\em the same as}
the restrictions to $\Z^n_+$ of the dual crystal operators
$\tilde E_i^*, \tilde F_i^*, \eps_i^*, \phi_i^*$ defined in 
\ref{crys}, so we will not need the latter notation again.

In fact, there are now so few possibilities that we can describe the
crystal graph explicitly. First suppose that $i = 0$.
Then the possible $i$-strings in the crystal graph are as follows:
\begin{itemize}
\item[(1)] $(\cdots)$;
\item[(2)]
$(\cdots,0^r,-1,\cdots) \stackrel{\tilde F_0}{\longrightarrow}
(\cdots,0^{r+1},\cdots) \stackrel{\tilde F_0}{\longrightarrow}
(\cdots,1,0^r,\cdots)$;
\item[(3)]
$(\cdots,1,0^r,-1,\cdots)$.
\end{itemize}
Here, $\cdots$ denotes some fixed entries different from $1,0,-1$
and $r \geq 0$.
Similarly for $i > 0$, the possible $i$-strings in the
crystal are as follows:
\begin{itemize}
\item[(1)] $(\cdots)$;
\item[(2)]
$(\cdots,i,\cdots) 
\stackrel{\tilde F_i}{\longrightarrow}
(\cdots,i+1,\cdots)$;
\item[(3)]
$(\cdots,-i-1,\cdots) 
\stackrel{\tilde F_i}{\longrightarrow}
(\cdots,-i,\cdots)$;
\item[(4)]
$(\cdots,i,\cdots,-i-1,\cdots)
\stackrel{\tilde F_i}{\longrightarrow}
(\cdots,i,\cdots,-i,\cdots)\\
\phantom{hello}\hspace{60mm}
\stackrel{\tilde F_i}{\longrightarrow}
(\cdots,i+1,\cdots,-i,\cdots)$;
\item[(5)]
$(\cdots,i+1,\cdots,-i-1,\cdots)$;
\item[(6)]
$(\cdots,-i,-i-1,\cdots)$;
\item[(7)]
$(\cdots,i+1,i,\cdots)$;
\item[(8)]
$(\cdots,i+1,i,\cdots,-i-1,\cdots) \stackrel{\tilde F_i}{\longrightarrow}
(\cdots,i+1,i,\cdots,-i,\cdots)$;
\item[(9)]
$(\cdots,i,\cdots,-i,-i-1,\cdots)
\stackrel{\tilde F_i}{\longrightarrow}
(\cdots,i+1,\cdots,-i,-i-1,\cdots)$;
\item[(10)] $(\cdots,i+1,i,\cdots,-i,-i-1,\cdots)$,
\end{itemize}
where again $\cdots$ denotes fixed entries different from $i,i+1,-i,-i-1$.
A crucial observation deduced from this analysis is
that {\em all} $i$-strings are of length $\leq 2$.

\begin{Lemma}\label{ccthm}
Let $\lambda \in \Z_+^{n}$ and $i \geq 0$.
\begin{itemize}
\item[(i)]
$E_i L_\lambda = [\eps_i(\la)]_i L_{\tilde E_i (\la)} + 
\sum_{\mu \in \Z_+^{n}}
w_{\mu,\la}^i L_\mu$
where $w_{\mu,\la}^i\in qq_i^{1 - \eps_i(\la)} \Z[q]$
is zero unless
$\phi_j(\mu) \leq \phi_j(\la)$ for all $j \geq 0$.
\item[(ii)]
$F_i L_\la = [\phi_i(\la)]_i L_{\tilde F_i (\la)} + 
\sum_{\mu \in \Z_+^{n}}
x_{\mu,\la}^i L_\mu$
where $x_{\mu,\la}^i\in qq_i^{1 - \phi_i(\la)} \Z[q]$
is zero unless
$\eps_j(\mu) \leq \eps_j(\la)$ for all $j \geq 0$.
\end{itemize}
\end{Lemma}

\begin{proof}
This is a special case of
Lemma~\ref{cthm}, since $\tilde E_i = \tilde E_i^*$,
$\tilde F_i = \tilde F_i^*, \dots$.
\end{proof}

\begin{Lemma}\label{thme}
Let $\lambda \in \Z_+^{n}$ and $i \geq 0$.
\begin{itemize}
\item[(i)]
$E_i U_\lambda = 
[\phi_i(\lambda)+1]_i U_{\tilde E_i(\lambda)} + 
\sum_{\mu \in \Z_+^{n}}
y_{\mu,\lambda}^i U_\mu$
where $y_{\mu,\lambda}^i\in qq_i^{1 - \phi_i(\mu)} \Z[q]$
is zero unless
$\eps_j(\mu) \geq \eps_j(\lambda)$ for all $j \geq 0$.
\item[(ii)]
$F_i U_\lambda = [\eps_i(\lambda)+1]_i U_{\tilde F_i(\lambda)} + 
\sum_{\mu \in \Z_+^{n}}
z_{\mu,\lambda}^i U_\mu$
where $z_{\mu,\lambda}^i\in qq_i^{1 - \eps_i(\mu)} \Z[q]$
is zero unless
$\phi_j(\mu) \geq \phi_j(\lambda)$ for all $j \geq 0$.
\end{itemize}
(In (i) resp. (ii), the first term on the right hand side 
should be omitted if
$\tilde E_i(\la)$ resp. 
$\tilde F_i(\la)$ is $\varnothing$.)
\end{Lemma}

\begin{proof}
Dualize Lemma~\ref{ccthm}.
\end{proof}

\begin{Corollary}\label{cor} Let $\la \in \Z^{n}_+$ and $i \geq 0$.
\begin{itemize}
\item[(i)]
If $\eps_i(\la) > 0$ then
$E_i U_\lambda = 
[\phi_i(\la)+1]_i U_{\tilde E_i (\la)}$.

\item[(ii)]
If $\phi_i(\la) > 0$ then
$F_i U_\lambda = 
[\eps_i(\la)+1]_i U_{\tilde F_i (\la)}$.
\end{itemize}
\end{Corollary}

\begin{proof}
We prove (i), (ii) being similar.
Lemma~\ref{thme} gives us that
$E_i U_\la = 
[\phi_i(\la)+1]_i U_{\tilde E_i(\la)} + 
\sum_{\mu\in \Z_+^{n}}
y_{\mu,\la}^i U_\mu$
where $y_{\mu,\la}^i$ 
belongs to $qq_i^{1 - \phi_i(\mu)} \Z[q]$
and is zero unless
$\eps_j(\mu) \geq \eps_j(\mu)$ for all $j\geq 0$.
Suppose that $y_{\mu,\la}^i \neq 0$ for some $\mu$.
By assumption, $\eps_i(\mu) \geq \eps_i(\la) \geq 1$,
so $\phi_i(\mu) \leq 1$ since all $i$-strings are of length $\leq 2$.
So $0 \neq y_{\mu,\la}^i \in q \Z[q]$.
But $y_{\mu,\la}^i$ is bar invariant, so this is a contradiction.
\end{proof}

\Point{\boldmath Computation of $U_\lambda$'s}\label{alg}
Now we explain a simple algorithm to compute $U_\lambda$.
Recall the definition of the degree of atypicality of $\lambda$ from
(\ref{atyp}). We will need the following:

\begin{Procedure}\label{alg2}\rm
Suppose we are given $\lambda \in \Z^n_+$ with $\#\lambda \geq 2$.
Compute $\mu \in \Z^n_+$ and an operator
$X_i
\in \{E_i, F_i\}_{i \geq 0}$ by following the instructions below starting
at step (0).
\begin{itemize}
\item[(0)\phantom{$'$}]
Choose the minimal $r \in \{1,\dots,n\}$ such that $\lambda_r+\lambda_s = 0$
for some $s > r$. Go to step (1).
\item[(1)\phantom{$'$}] 
If $r > 1$ and $\la_r = \lambda_{r-1} - 1$, replace $r$ by $(r-1)$
and repeat step (1). Otherwise, go to step (2).
\item[(2)\phantom{$'$}] If $\lambda_r + \la_s + 1 = 0$ for some (necessarily unique)
$s \in \{1,\dots,n\}$
go to step (1)$'$. 
Otherwise, set $X_i = E_{\lambda_r}$ and 
$\mu = \lambda+ \delta_r$. Stop.
\item[(1)$'$]
If $s < n$ and $\la_s= \lambda_{s+1} + 1$, replace $s$ by $(s+1)$
and repeat step (1)$'$. Otherwise, go to step (2)$'$.
\item[(2)$'$] If $\la_r + \lambda_s - 1 = 0$ for some (necessarily unique)
$r \in \{1,\dots,n\}$ 
go to step (1). 
Otherwise, set $X_i = F_{-\la_s}$ and 
$\mu = \lambda - \delta_s$. Stop.
\end{itemize}
\end{Procedure}

The following lemma follows immediately from the nature of the
above procedure and Corollary~\ref{cor}.

\begin{Lemma}\label{props}
Take $\lambda \in \Z^n_+$
with $\#\lambda \geq 2$.
Apply Procedure~\ref{alg2} to get
$\mu \in \Z^n_+$ and $X_i \in \{E_i, F_i\}_{i \geq 0}$.
Then, $\#\mu \leq \#\la$ and $X_i U_\mu = U_\la$.
Moreover, after at most $(n-1)$ repetitions of the procedure,
the atypicality must get strictly smaller.
Hence after finitely many recursions, the 
procedure reduces $\lambda$ to a {typical} weight.
\end{Lemma}

Lemma~\ref{props}
implies the following algorithm for computing $U_\lambda$.
If $\lambda$ is typical then $U_\lambda = F_\lambda$,
since such $\lambda$'s are maximal in the Bruhat ordering.
Otherwise, apply Procedure~\ref{alg2} to get $\mu \in \Z^n_+$ and $X_i
\in \{E_i, F_i\}_{i \geq 0}$.
Since the procedure always reduces $\lambda$ to a typical weight in 
finitely many steps, we may assume $U_\mu$ is known recursively. Then
$U_\lambda = X_i U_\mu$.

\begin{Example}\label{longone}\rm
Let us compute $U_{(5,3,2,1,0,0,-1,-4,-6)}$. Apply
Procedure~\ref{alg2} repeatedly to get the following sequence of weights:
\begin{align*}
(5,3,2,1,0,0,-1,-4,-7),
\qquad&(6,3,2,1,0,0,-1,-4,-7),\\
(6,3,2,1,0,0,-1,-5,-7),
\qquad&(6,4,2,1,0,0,-1,-5,-7), \\
(6,4,3,1,0,0,-1,-5,-7),
\qquad&(6,4,3,2,0,0,-1,-5,-8),\\ 
(7,4,3,2,0,0,-1,-6,-8),
\qquad&(7,5,3,2,0,0,-1,-6,-8),\\
(7,5,4,2,0,0,-1,-6,-8),
\qquad&(7,5,4,3,0,0,-1,-6,-8), \\
(7,5,4,3,0,0,-2,-6,-8),
\qquad&(7,5,4,3,1,0,-2,-6,-8).
\end{align*}
Hence,
\begin{align*}
U_{(5,3,2,1,0,0,-1,-4,-6)}
&=
F_6E_5F_4E_3E_2E_1F_7E_6F_5E_4E_3E_2F_1E_0 F_{(7,5,4,3,1,0,-2,-6,-8)}\\
&= 
F_{(5,3,2,1,0,0,-1,-4,-6)}
+
q^2F_{(7,5,3,2,0,0,-4,-6,-7)}\\
&\qquad+
(q+q^3)F_{(8,5,3,2,1,-1,-4,-6,-8)}\\
&\qquad+
(q^3+q^5)F_{(8,7,5,3,2,-4,-6,-7,-8)}.
\end{align*}
\end{Example}

\begin{Example}\rm
Using the algorithm, one computes for $n = 3$ that:
\begin{align*}
U_{(a,b,c)} &= F_{(a,b,c)}&(a+b,a+c,b+c\neq 0)\\
U_{(a,b,-b)} &= F_{(a,b,-b)} + q^2 F_{(a,b+1,-b-1)}&(a \geq b+2, b \neq 0)\\
U_{(b+1,b,-b)} &= F_{(b+1,b,-b)} + q^2 F_{(b+2,b+1,-b-2)}&(b \neq 0)\\
U_{(a,-a,-b)} &= F_{(a,-a,-b)} + q^2 F_{(a+1,-a-1,-b)}&(b \geq a+2, a \neq 0)\\
U_{(a,-a,-a-1)} &= F_{(a,-a,-a-1)} + q^2 F_{(a+2,-a-1,-a-2)}&(a \neq 0)\\
U_{(a,b,-a)} &= F_{(a,b,-a)} + q^2 F_{(a+1,b,-a-1)}&(a > b > -a)\\
U_{(a,0,0)} &= F_{(a,0,0)} + (q+q^3) F_{(a,1,-1)}&(a \geq 2)\\
U_{(1,0,0)} &= F_{(1,0,0)} + (q+q^3) F_{(2,1,-2)}&\\
U_{(0,0,-b)} &= F_{(0,0,-b)} + (q+q^3) F_{(1,-1,-b)}&(b \geq 2)\\
U_{(0,0,-1)} &= F_{(0,0,-1)} + (q+q^3) F_{(2,-1,-2)}&\\
U_{(0,0,0)} &= F_{(0,0,0)}+(q-q^5) F_{(1,0,-1)} + (q+q^3) F_{(2,0,-2)}\!\!\!\!\!\!\!\!\!\!
\end{align*}
\end{Example}

\Point{Specialization}\label{zforms}
When we specialize at $q = 1$, various things become simpler to compute.
To formalize the specialization process, we will work with the
$\Z[q,q^{-1}]$-lattices
\begin{align*}
\E^n_{\Z[q,q^{-1}]} &= \sum_{\lambda \in \Z_+^n} \Z[q,q^{-1}] L_\lambda
= \sum_{\lambda \in \Z_+^n} \Z[q,q^{-1}] E_\lambda,\\
\F^n_{\Z[q,q^{-1}]} &= \sum_{\lambda \in \Z_+^n} \Z[q,q^{-1}] U_\lambda
\subset \sum_{\lambda \in \Z^n_+} \Z[q,q^{-1}] F_\lambda.
\end{align*}
Also let $\FF^n_{\Z[q,q^{-1}]}$ denote the
completion of $\F^n_{\Z[q,q^{-1}]}$, i.e. its closure in $\FF^n$.
We need this because the element $F_\lambda \in \F^n$ need not belong to
the lattice $\F^n_{\Z[q,q^{-1}]}$, though it always
belongs to $\FF^n_{\Z[q,q^{-1}]}$. Indeed, the
$\{F_\lambda\}_{\lambda \in \Z^n_+}$ form a topological basis 
for $\FF^n_{\Z[q,q^{-1}]}$.

\begin{Lemma}\label{argb}
$\E^n_{\Z[q,q^{-1}]}, \F^n_{\Z[q,q^{-1}]}$ and $\FF^n_{\Z[q,q^{-1}]}$ 
are modules
over $\U_{\Z[q,q^{-1}]}$.
\end{Lemma}

\begin{proof}
The $\Z[q,q^{-1}]$-lattice in
$\TT^n$ generated by $\{M_\lambda\}_{\lambda \in \Z^n}$
is invariant under
$\U_{\Z[q,q^{-1}]}$.
Combining this with Corollary~\ref{fne} and Lemma~\ref{ems}, 
we deduce that each
$E_i^{(r)} E_\lambda$ can be expressed as 
a (possibly infinite) 
$\Z[q,q^{-1}]$-linear combination of $E_\mu$'s. 
Actually, it is a finite linear combination by Lemma~\ref{expact}.
Similarly, each $F_i^{(r)} E_\lambda$ is a finite
$\Z[q,q^{-1}]$-linear combination of $E_\mu$'s. So
$\E^n_{\Z[q,q^{-1}]}$ is a $\U_{\Z[q,q^{-1}]}$-module.

Considering instead the
$\Z[q,q^{-1}]$-lattice in
$\TT^n$ generated by $\{N_\lambda\}_{\lambda \in \Z^n}$ and passing
to the quotient $\FF^n$, we also get easily
that each $E_i^{(r)} F_\la$ and each
$F_i^{(r)} F_\la$ is a finite
$\Z[q,q^{-1}]$-linear combination of $F_\mu$'s,
hence $\FF^n_{\Z[q,q^{-1}]}$ is a $\U_{\Z[q,q^{-1}]}$-module.
This means in particular that each
$E_i^{(r)} U_\la$ and each
$F_i^{(r)} U_\la$ is a (possibly infinite)
$\Z[q,q^{-1}]$-linear combination of $U_\mu$'s.
Finally to show that $\F^n_{\Z[q,q^{-1}]}$ is a 
$\U_{\Z[q,q^{-1}]}$-module we need to show 
$E_i^{(r)} U_\la$ and
$F_i^{(r)} U_\la$ are actually finite
$\Z[q,q^{-1}]$-linear combination of $U_\mu$'s.
We explain the argument for $E_i^{(r)} U_\la$ only.

By the algorithm described in the previous subsection, 
$U_\la$ is a finite linear combination of $F_\mu$'s, hence
$E_i^{(r)} U_\la$ is a bar invariant, 
finite linear combination of $F_\mu$'s. Say 
$$
E_i^{(r)} U_\la = \sum_{\mu \in \Z^n} f_\mu(q) F_\mu
$$
for polynomials $f_\mu(q) \in \Z[q,q^{-1}]$. 
Let $d$ be minimal such that all $f_\mu(q)$ belong to $q^{-d} \Z[q]$.
If $d < 0$ then the right hand side is
a $q \Z[q]$-linear combination of $F_\mu$'s, hence since it also
a bar invariant combination of $U_\mu$'s it must be zero.
Otherwise, 
define $g_\mu(q)$ to be the unique bar invariant polynomial with
$g_\mu(q) \equiv f_\mu(q) \pmod{q\Z[q]}$
and subtract $\sum_\mu g_\mu(q) U_\mu$ from the right hand side,
summing over the finitely many $\mu$ such that $\deg g_\mu(q) = d$.
The result is a bar invariant, finite 
$q^{1-d} \Z[q]$-linear combination of $F_\mu$'s.
Now repeat the process to reduce the expression to zero after finitely many 
steps.
\end{proof}

Define
\begin{align*}
\E^n_\Z &= \Z \otimes_{\Z[q,q^{-1}]} \E^n_{\Z[q,q^{-1}]},\\
\F^n_\Z &= \Z \otimes_{\Z[q,q^{-1}]} \F^n_{\Z[q,q^{-1}]},\\
\FF^n_\Z &= \Z \otimes_{\Z[q,q^{-1}]} \FF^n_{\Z[q,q^{-1}]},
\end{align*}
where we are viewing $\Z$ as a $\Z[q,q^{-1}]$-module so that $q$ acts as $1$.
For $\lambda \in \Z_+^n$, we 
will denote the elements 
$1 \otimes E_\lambda \in \E^n_{\Z}, 1 \otimes L_\lambda \in \E^n_{\Z}$,
$1 \otimes F_\lambda \in \FF^n_{\Z}$ and $1 \otimes U_\lambda \in \F^n_{\Z}$
by 
$E_\lambda(1)$, $L_\lambda(1)$, $F_\lambda(1)$ and $U_\lambda(1)$,
respectively.
By Lemma~\ref{argb}, $\E^n_{\Z}, \F^n_{\Z}$ and $\FF^n_{\Z}$ are modules 
over $\U_{\Z} := \Z \otimes_{\Z[q,q^{-1}]} \U_{\Z[q,q^{-1}]}$.
In their action on these lattices, the elements
\begin{equation}\label{generators}
E_i^{(r)} = 1 \otimes E_i^{(r)},\quad
F_i^{(r)} = 1 \otimes F_i^{(r)}
,\quad
\binom{H_i}{r} = 1 \otimes \left[
\begin{array}{c}
K_i\\r\end{array}\right]
\end{equation}
of $\U_{\Z}$ satisfy the defining relations of the usual generators
of the Kostant $\Z$-form
for the universal enveloping algebra of type $\mathfrak{b}_{\infty}$.

We note that at $q = 1$, the
defining relations for $\F_{\Z}^n$ simplify to the following:
\begin{align*}
v_a \wedge v_a&= 0&(a \neq 0)\\
v_a \wedge v_b &= - v_b \wedge v_a&(a > b, a+b\neq 0)\\
v_{a}\wedge v_{-a} &= -v_{-a}\wedge v_a + (-1)^a v_0\wedge v_0
&(a \geq 1)
\end{align*}
The action of $\U_{\Z}$ on $\FF^n_{\Z}$ is given explicitly by
the formulae
\begin{align}\label{a1}
E_i F_\lambda(1) &= \sum_{\substack{r\:\text{with}\:\la - \delta_r \in \Z^n_+ \\ \sigma_r = -, -+\:\text{or}\:-- }}
b_{\la,r} F_{\la-\delta_r}(1),\\
F_i F_\lambda(1) &= \sum_{\substack{r\:\text{with}\:\la + \delta_r \in \Z^n_+ \\ \sigma_r = +, -+\:\text{or}\:++ }}
b_{\la,r} F_{\la+\delta_r}(1),\label{a2}
\end{align}
where $(\sigma_1,\dots,\sigma_n)$ is the $i$-signature
of $\la\in\Z^n_+$ defined according to (\ref{isig})
and 
$$
b_{\la,r} = \left\{\begin{array}{ll}
1+(-1)^{z(\la)}&\hbox{if $\sigma_r = -+$,}\\
1&\hbox{otherwise.}
\end{array}
\right.
$$

The pairing $\langle.,.\rangle$ between $\EE^n$ and 
$\FF^n$ defined earlier induces a pairing at $q = 1$ between
$\E^n_{\Z}$ and $\FF^n_{\Z}$ with
\begin{equation}
\langle L_{-w_0\la}(1),U_\mu(1) \rangle = \delta_{\la,\mu}
\end{equation}
for all $\la,\mu \in \Z^n_+$.
Alternatively, as follows immediately from (\ref{ude}) and (\ref{edef})
taken at $q = 1$, we have that
\begin{equation}\label{formy}
\langle E_{-w_0\la}(1), F_\mu(1) \rangle = \delta_{\la,\mu}
\end{equation}
for all $\la,\mu \in \Z^n_+$.

The final theorem of the section gives an explicit formula for the
coefficients $u_{\mu,\lambda}(q)$ at $q = 1$:

\begin{Theorem}\label{mt}
Let $\lambda \in \Z^n_+$ and $p := \lfloor \#\lambda / 2 \rfloor$.
Choose $1\leq r_1 < \dots < r_p < s_p < \dots < s_1 \leq n$ 
such that $\lambda_{r_q} + \lambda_{s_q} = 0$ for all
$q = 1,\dots,p$.
Let $I_0 =\{|\la_1|,\dots,|\la_n|\}$.
For $q = 1,\dots,p$,
define $I_q$ and $k_q$
inductively according to the following rules:
\begin{itemize}
\item[(1)]
if $\la_{r_q} > 0$, let $k_q$ be the smallest
positive integer with $\la_{r_q} + k_q \notin I_{q-1}$, and set
$I_q = I_{q-1} \cup \{\la_{r_q} + k_q\}$;
\item[(2)]
if $\la_{r_q} = 0$, let $k_q$ and $k_q'$ be the smallest positive
integers with $k_q, k_q' \notin I_{q-1}$,
$k_q < k_q'$ if $z(\la)$ is even and
$k_q > k_q'$ if $z(\la)$ is odd, and set
$I_q = I_{q-1} \cup \{k_q,k_q'\}$.
\end{itemize}
Finally, for each $\theta = (\theta_1,\dots,\theta_p) \in \{0,1\}^p$,
let $\RR_{\theta}(\lambda)$ denote the unique element of $\Z^n_+$
lying in the same $S_n$-orbit as the weight
$\la + \sum_{q=1}^p \theta_q k_q (\delta_{r_q} - \delta_{s_q})$.
Then, 
$$
U_\lambda(1) = \sum_{\theta} 2^{(z(\lambda)-z(\RR_{\theta}(\lambda)))/2}
F_{\RR_{\theta}(\lambda)}(1),
$$
summing over all 
$\theta = (\theta_1,\dots,\theta_p)\in\{0,1\}^p$.
\end{Theorem}

\begin{proof}
Use the algorithm explained in \ref{alg}, (\ref{a1})--(\ref{a2}) 
and induction. 
See \cite[Theorem 3.34(i)]{KL} for a similar argument.
\end{proof}

\begin{Corollary} For all $\lambda \in \Z^n_+$,
$\sum_{\mu \in \Z^n_+} u_{\mu,\lambda}(1) = 
2^{(\#\lambda - z(\lambda))/2} 3^{\lfloor z(\lambda)/2\rfloor}.$
\end{Corollary}

\begin{Example}\rm
Let us compute $U_{(5,3,2,1,0,0,-1,-4,-6)}(1)$ using the theorem,
recall Example~\ref{longone}.
We have that $p=2$ and $(i_1, i_2, j_2, j_1) = (4,5,6,7)$,
then get $k_1 = 6$ (hence the entries $(1,-1)$ change to
$(7,-7)$ and $(k_2,k_2') = (8,9)$ (hence the entries $(0,0)$ change
to $(8,-8)$). Hence
\begin{multline*}
U_{(5,3,2,1,0,0,-1,-4,-6)}(1) = 
F_{(5,3,2,1,0,0,-1,-4,-6)}(1)
+
F_{(7,5,3,2,0,0,-4,-6,-7)}(1)\\
+
2F_{(8,5,3,2,1,-1,-4,-6,-8)}(1)
+
2F_{(8,7,5,3,2,-4,-6,-7,-8)}(1).
\end{multline*}
For comparison,
\begin{multline*}
U_{(5,3,2,1,0,0,0,-1,-4,-6)}(1) = 
F_{(5,3,2,1,0,0,0,-1,-4,-6)}(1)
+
F_{(7,5,3,2,0,0,0,-4,-6,-7)}(1)\\
+
2F_{(9,5,3,2,1,0,-1,-4,-6,-9)}(1)
+
2F_{(9,7,5,3,2,0,-4,-6,-7,-9)}(1).
\end{multline*}
\end{Example}

\begin{Example}\rm We have by the theorem that:
\begin{align*}
U_{(0,0)}(1) &= F_{(0,0)}(1) + 2 F_{(1,-1)}(1),\\
U_{(0,0,0)}(1) &= F_{(0,0,0)}(1) + 2 F_{(2,0,-2)}(1),\\
U_{(0,0,0,0)}(1) &= F_{(0,0,0,0)}(1) + 2 F_{(1,0,0,-1)}(1)
 + 2 F_{(3,0,0,-3)}(1) + 4 F_{(3,1,-1,-3)}(1),\\
U_{(1,-1)}(1) &= F_{(1,-1)}(1) + F_{(2,-2)}(1),\\
U_{(1,0,-1)}(1) &= F_{(1,0,-1)}(1) + F_{(2,0,-2)}(1),\\
U_{(1,0,0,-1)}(1) &= F_{(1,0,0,-1)}(1) + F_{(2,0,0,-2)}(1)
+ 2F_{(3,1,-1,-3)}(1) + 2 F_{(3,2,-2,-3)}(1).\\\intertext{Hence by (\ref{edef}),}
E_{(1,-1)}(1) &= L_{(1,-1)}(1) + 2 L_{(0,0)}(1),\\
E_{(1,0,-1)}(1) &= L_{(1,0,-1)}(1),\\
E_{(1,0,0,-1)}(1) &= L_{(1,0,0,-1)}(1) + 2 L_{(0,0,0,0)}(1),\\
E_{(2,-2)}(1) &= L_{(2,-2)}(1),\\
E_{(2,0,-2)}(1) &= L_{(2,0,-2)}(1) + L_{(1,0,-1)}(1)+ 2L_{(0,0,0)}(1),\\
E_{(2,0,0,-2)}(1) &= L_{(2,0,0,-2)}(1) + L_{(1,0,0,-1)}(1).
\end{align*}
\end{Example}
\vspace{1mm}

\section{Character formulae}\label{cf}

\Point{\boldmath Representations of $Q(n)$}\label{reps}
Now we are ready to introduce the supergroup $G = Q(n)$ into the picture.
For the remainder of the article, we will work
over the ground field $\C$, and all objects (superalgebras, 
superschemes, \dots)
will be defined over $\C$ without further mention.
We refer the reader to \cite{BKI, BII} for a fuller account of
the basic results concerning the representation theory of $G$ summarized 
here,
most of which were proved originally by Penkov in \cite{penkov}.

By definition \cite[$\S$3]{BKI}, $G$ is the functor from 
the category of commutative superalgebras to the category of groups
defined on a superalgebra $A$
so that $G(A)$ is the group of all invertible $2n \times 2n$ matrices
of the form
\begin{equation}\label{form}
g = \left(
\begin{array}{r|r}
S\phantom{'}&S'\\\hline
-S'&S\phantom{'}
\end{array}
\right)
\end{equation}
where $S$ is an $n \times n$ matrix with entries in $A_{\0}$
and $S'$ is an $n \times n$ matrix with entries in $A_{\1}$.
The {underlying purely even group} $G_{\ev}$ of $G$ is by definition
the functor from superalgebras to groups 
with $G_\ev(A) := G(A_{\0})$. In our case, $G_\ev(A)$ consists
of all matrices in $G(A)$ of the form (\ref{form}) with $S'=0$,
so $G_\ev \cong GL(n)$.
We also need the {\em Cartan subgroup} $H$ of $G$ 
defined on a commutative superalgebra $A$ so that $H(A)$ consists
of all matrices in $G(A)$ with $S, S'$ being
diagonal matrices, and the {\em negative Borel subgroup} $B$ of $G$
defined so that $B(A)$ consists of all
matrices in $G(A)$ with $S, S'$ being lower  triangular. 

Let $T = H_{\ev}$ 
be the usual maximal torus of $G_\ev$ consisting of diagonal matrices.
The character group 
$X(T) = \operatorname{Hom}(T, \mathbb{G}_m)$ is the free abelian
group on generators $\delta_1,\dots,\delta_n$, where $\delta_r$
picks out the $r$th diagonal entry of a diagonal matrix.
We will always identify $X(T)$ with $\Z^n$,
the tuple $\lambda = (\lambda_1,\dots,\lambda_n) \in \Z^n$
corresponding to the character $\sum_{r=1}^n \lambda_r\delta_r \in X(T)$.
The root system associated to $G_{\ev}$ is denoted  
$R = R^+ \cup (-R^+)$, where
$R^+ = \{\delta_r - \delta_s\:|\:1 \leq r< s \leq n\} \subset \Z^n$.
The dominance ordering on $\Z^n$ is defined by 
$\la \leq \mu$ if and only if $\mu - \la$ is a sum of positive roots.
This is not the same as the Bruhat ordering $\preceq$
defined in \ref{bo}, though we have by Lemma~\ref{bruhat}
that $\la \preceq \mu$ implies $\la \leq \mu$.

A representation of $G$ means 
a natural transformation $\rho:G \rightarrow GL(M)$
for some vector superspace $M$, where $GL(M)$ is the supergroup
with $GL(M, A)$ being equal to the group of all even automorphisms of 
the $A$-supermodule $M \otimes A$, for each commutative superalgebra $A$.
Equivalently, as with group schemes \cite[I.2.8]{Jantzen}, $M$ is  
a right $k[G]$-comodule, i.e. there is an even
structure map $\eta: M \to M \otimes k[G]$ satisfying the usual comodule 
axioms. 
We will usually refer to such an $M$ as a {\em $G$-supermodule}.
Let $\mathcal C_n$ denote the category of all finite dimensional
$G$-supermodules.
Note we allow arbitrary (not necessarily homogeneous) morphisms
in $\mathcal C_n$. We write $\Pi$ for the parity change functor, 
and denote the dual of a finite dimensional $G$-supermodule
$M$ by $M^*$. There is another natural duality $M \mapsto M^\tau$
on $\mathcal C_n$, 
see e.g. \cite[$\S$10]{BKI} for the definition.

There are two sorts of irreducible $G$-supermodule:
either {\em type $\Mtype$} if $\End_G(M)$ is one dimensional or 
{\em type $\Qtype$} if $\End_G(M)$ is two dimensional.
For example, we have the \emph{natural representation} $V$ of $G$, 
the vector superspace on basis
$v_1,\dots,v_n, v_1',\dots,v_n'$, where $v_r$ is even and $v_r'$ is odd.
For a superalgebra $A$, we identify elements of $V \otimes A$ with column vectors
$$
\sum_{r=1}^n (v_r \otimes a_r + v_r' \otimes a_r')
\longleftrightarrow
\left(
\begin{array}{l}
a_1\\\vdots\\a_n\\a_1'\\\vdots\\a_n'\end{array}
\right).
$$
Then, the action of $G(A)$ on $V \otimes A$ defining the supermodule structure
is the obvious action on column vectors
by left multiplication.
The map \begin{equation}\label{jdef}
J:V \rightarrow V,\qquad
v_r \mapsto v_r', v_r'\mapsto -v_r
\end{equation}
is an odd automorphism of $V$ as a $G$-supermodule.
Hence, $V$ is irreducible of type $\Qtype$.

For $\lambda \in \Z^n$, 
we write $x^\lambda = x_1^{\lambda_1} \dots x_n^{\lambda_n}
\in \Z[x_1^{\pm 1}, \dots, x_n^{\pm 1}]$.
Given an $H$-supermodule $M$,
we let $M_\lambda$ denote its $\lambda$-weight space with respect to 
the torus $T$.
Identifying the Weyl group associated to $G_\ev$ with the symmetric group
$S_n$,
the character 
$$
\ch M := \sum_{\lambda \in \Z^n} (\dim M_\lambda) x^\lambda
$$ 
of a finite dimensional $G$-supermodule $M$
is naturally $S_n$-invariant, so is an element of the ring
$\Z[x_1^{\pm 1},\dots,x_n^{\pm 1}]^{S_n}$
of symmetric functions.

For every $\lambda \in \Z^n$, there is by \cite[Lemma 6.4]{BKI} 
a unique irreducible $H$-supermodule denoted $\u(\lambda)$ with character
$2^{\lfloor (h(\lambda)+1) / 2 \rfloor} x^{\lambda}$,
where $h(\lambda) = n - z(\la)$ denotes the number of $r = 1,\dots,n$
for which $\lambda_r \neq 0$.
It is an irreducible $H$-supermodule of type $\Mtype$ if $h(\lambda)$ is
even, type $\Qtype$ otherwise.
We will often
regard $\u(\lambda)$ instead as an irreducible $B$-supermodule via
the obvious epimorphism $B \rightarrow H$.
Introduce the induced supermodule
\begin{equation}\label{h0def}
H^0(\lambda) := \ind_B^G \u(\lambda) = H^0 (G / B, \mathscr L(\u(\lambda))),
\end{equation}
see \cite[$\S$6]{BKI} and \cite[$\S$2]{BII} for the detailed construction.
The following theorem 
is due to Penkov \cite[Theorem 4]{penkov},
see also \cite[Theorem 6.11]{BKI}.

\begin{Theorem}\label{class}
For $\lambda \in\Z^n$, $H^0(\lambda)$ is non-zero
if and only if $\lambda \in \Z^n_+$. In that
case $H^0(\lambda)$ has a unique irreducible submodule denoted $L(\lambda)$.
Moreover, 
\begin{itemize}
\item[(i)] $\{L(\lambda)\}_{\lambda \in \Z^n_+}$ is a complete set of pairwise
non-isomorphic irreducible $G$-supermodules; 
\item[(ii)]
$L(\lambda)$ is of the same type as $\u(\lambda)$,
i.e. type $\Mtype$ if $h(\lambda)$ is even, type $\Qtype$ otherwise;
\item[(iii)] 
$\ch L(\lambda) = 2^{\lfloor (h(\lambda)+1)/2 \rfloor} x^\lambda + (*)$
where $(*)$ is a linear combination of $x^\mu$ for $\mu < \lambda$;
\item[(iv)] $L(\lambda)^* \cong L(-w_0\lambda)$
and $L(\lambda)^\tau \cong L(\la)$.
\end{itemize}
\end{Theorem}

\Point{Euler characteristics}\label{ecsec}
Let $K(\mathcal C_n)$ denote the Grothendieck group of the 
superadditive category $\mathcal C_n$, in the sense of \cite[$\S$2-c]{BK}.
So by Theorem~\ref{class}(i), $K(\mathcal C_n)$ 
is the free abelian group on basis
$\left\{[L(\lambda)]\right\}_{\lambda \in \Z^n_+}$.
Since $\ch$ is additive on short exact sequences, there is an induced
map
$$
\ch:K(\mathcal C_n) 
\rightarrow \Z[x_1^{\pm 1},\dots,x_n^{\pm 1}]^{S_n},\qquad
[M] \mapsto \ch M.
$$
Theorem~\ref{class}(iii) implies that the irreducible characters
are linearly independent, hence this map is injective.
By Theorem~\ref{class}(iv), 
the duality $*$ on finite dimensional $G$-supermodules
induces an involution $*:K(\mathcal C_n) \rightarrow K(\mathcal C_n)$
with $[L(\lambda)]^* = [L(-w_0\lambda)]$.
On the other hand, the duality
$\tau$ leaves characters invariant, so gives the identity
map at the level of Grothendieck groups.

For each $\lambda \in \Z^n$, we have the higher cohomology supermodules
\begin{equation}
H^i(\lambda) := R^i \ind_P^G \u(\lambda)
=
H^i (G / P, \mathscr L(\u(\lambda)))
\end{equation}
where $P$ is the largest parabolic subgroup of $G$ to which the $B$-supermodule
$\u(\lambda)$ can be lifted, 
we refer to \cite[$\S$4]{BII} for details.
In particular, it is known that each $H^i(\lambda)$
is finite dimensional and is zero for $i \gg 0$, so the Euler characteristic
\begin{equation}\label{euler}
[E(\lambda)] := \sum_{i \geq 0} (-1)^i [H^i(\lambda)]
\end{equation}
is a well-defined element of the Grothendieck group $K(\mathcal C_n)$.
Its character
\begin{equation}\label{eulerch}
\ch E(\lambda) := \sum_{i \geq 0} (-1)^i \ch H^i(\lambda)
\end{equation}
can be computed explicitly by a method going back at least to Penkov
\cite[$\S$2.3]{p1}.
We need to recall the definition of
{\em Schur's $P$-function} $\P_\lambda$
for $\la \in \Z^n_+$:
\begin{equation}\label{pdef}
\P_\la = 
\sum_{w \in S_n / S_\lambda}
w \Bigg( x^\la \prod_{\substack{1 \leq i < j \leq n\\\la_i > \la_j}} 
\frac{1+ x_i^{-1}x_j}{1 - x_i^{-1}x_j}\Bigg),
\end{equation}
where $S_n / S_\lambda$ denotes the set of minimal length coset
representatives for the stabilizer $S_\lambda$ of $\lambda$ in $S_n$.
Note that $\P_\la$ is equal to the Hall-Littlewood symmetric function
\begin{equation}\label{hl}
\P_\la(t) = 
\sum_{w \in S_n / S_\lambda}
w \Bigg( x^\la \prod_{\substack{1 \leq i < j \leq n\\\la_i > \la_j}} 
\frac{x_i - tx_j}{x_i - x_j}\Bigg),
\end{equation}
evaluated at $t = -1$, see \cite[III(2.2)]{Mac} 
and \cite[III.8]{Mac}
(actually, Macdonald only describes the case when all $\la_i 
\geq 0$, but everything easily extends to $\la_i \in \Z$).
We note the following combinatorial lemma.

\begin{Lemma}\label{prodform}
$(x_1+\dots+x_n)p_\lambda  =
\sum_r p_{\lambda+\delta_r}$,
where the sum is over all $r=1,\dots,n$ such that
$\lambda+\delta_r \in \Z^n_+$ and 
moreover $\la_r \neq -1$ if $z(\la)$ is odd.
\end{Lemma}

\begin{proof}
Let $\la \in \Z^n$ with $\la_1\geq\dots\geq \la_n$.
Following the proof of \cite[III(3.2)]{Mac}, one shows that
$$
(x_1+\dots+x_n)p_\la(t)  = \sum_{r}
(1+t+\dots+t^{n_r})
 p_{\la+\delta_r}(t),
$$
where the sum is over all $r=1,\dots,n$ with
$\la_r < \la_{r-1}$ if $r > 1$, and
$n_r$ denotes $\#\{s\:|\:1\leq s\leq n, \la_s=\la_r+1\}$.
The lemma follows on setting $t = -1$.
\end{proof}

The following theorem is \cite[Proposition 1]{PS}, see \cite[Theorem 4.7]{BII}
for a more detailed exposition.

\begin{Theorem}\label{ecdef} For each $\la \in \Z^n_+$,
$\ch E(\la) = 
2^{\lfloor{(h(\lambda)+1)/2}\rfloor} \P_\la$.
\end{Theorem} 

This shows in particular that
$\ch E(\la)$
equals $2^{\lfloor{(h(\lambda)+1)/2}\rfloor}x^\la + (*)$ where $(*)$
is a linear combination of $x^\mu$ for $\mu < \la$.
Comparing with Theorem~\ref{class}(iii), we deduce that for
each $\la \in \Z^n_+$,  the {\em decomposition numbers} $d_{\mu,\la}$
defined from
\begin{equation}\label{split}
[E(\mu)]
=
\sum_{\la \in \Z^n_+} d_{\mu,\la} [L(\la)]
\end{equation}
satisfy $d_{\mu,\mu} = 1$ and $d_{\mu,\la} = 0$ 
for $\la\not\leq\mu$. Hence,
$\left\{[E(\la)]\right\}_{\la \in \Z^n_+}$
gives another natural basis for the Grothendieck group $K(\mathcal C_n)$.
Also define the {\em inverse decomposition numbers} $d_{\mu,\la}^{-1}$ from
\begin{equation}\label{split2}
[L(\mu)]
=
\sum_{\la \in \Z^n_+} d^{-1}_{\mu,\la} [E(\la)].
\end{equation}
Again we have that $d^{-1}_{\mu,\mu} = 1$ and $d^{-1}_{\mu,\la} = 0$ 
for $\la\not\leq\mu$.

\Point{\boldmath Category $\mathcal O_n$}
The Lie superalgebra $\mathfrak{g}$ of $G$ can be identified with 
the Lie superalgebra 
$\mathfrak{q}(n)$ of all
matrices of the form
\begin{equation}\label{form2}
x = \left(
\begin{array}{l|l}S&S'\\\hline S'&S\end{array}
\right)
\end{equation}
under the supercommutator $[.,.]$,
where $S$ and $S'$ are $n \times n$ matrices over $\C$ and
such a matrix is even if $S' = 0$ or odd if $S = 0$, 
see \cite[$\S$4]{BKI}.
We will let $e_{r,s}$ resp. $e_{r,s}'$ denote the even resp.
odd matrix unit, i.e. the matrix of the form (\ref{form2}) 
with the $rs$-entry of $S$ resp. $S'$ equal to $1$ and all other entries
equal to zero.
We will abbreviate $h_r := e_{r,r}, h_r' := e_{r,r}'$.
Let $\mathfrak h$ be the {\em Cartan subalgebra} of $\mathfrak g$
spanned by $\{h_r, h_r'\:|\:1 \leq r \leq n\}$,
and let $\mathfrak b$ be the {\em positive Borel
subalgebra} spanned by $\{e_{r,s}, e_{r,s}'\:|\:1 \leq r \leq s \leq n\}$.
Note $\mathfrak h$ is the Lie superalgebra of the
subgroup $H$, but perversely $\mathfrak b$ is not the Lie superalgebra
of $B$ since that consisted of lower triangular matrices!

For any $\lambda \in \Z^n$ and an $\mathfrak h$-supermodule $M$,
we define the $\lambda$-weight space $M_\lambda$ of $M$ by
$$
M_\lambda = \{m \in M\:|\:h_i m = \lambda_i m\hbox{ for each $i=1,\dots,n$}\}.
$$
When $M$ is an $H$-supermodule viewed as an $\mathfrak h$-supermodule
in the canonical way, the notion of weight space defined here
agrees with the earlier one.
Let $\mathcal O_n$ denote the category of all finitely generated
$\mathfrak g$-supermodules 
$M$ that are locally finite dimensional as $\mathfrak b$-supermodules and satisfy
$$
M = \bigoplus_{\lambda \in \Z^n} M_\lambda,
$$
cf. \cite{BGG}.
By \cite[Corollary 5.7]{BKI}, we can identify the category
$\mathcal C_n$ with the full subcategory of $\mathcal O_n$
consisting of all the finite dimensional objects.

\begin{Lemma}\label{filt}
Every $M \in \mathcal O_n$ has a composition series.
\end{Lemma}

\begin{proof}
Since the universal enveloping superalgebra
$U(\mathfrak g)$ of $\mathfrak g$ is Noetherian, every finitely generated
$\mathfrak g$-supermodule $M$ admits a descending
filtration $M = M_0 > M_1 > \dots$ with 
each $M_{i-1}/M_i$ simple.
We need to show this filtration is of finite length
in case $M \in \mathcal O_n$.
It suffices for this to show that the restriction of $M$ 
to $\mathfrak g_{\0} \cong \mathfrak{gl}(n)$ has a composition series.
Since $U(\mathfrak g)$ is a free left $U(\mathfrak g_{\0})$-supermodule
of finite rank, $M$ is still finitely generated when viewed as
a $\mathfrak g_{\0}$-module, and clearly it is
locally finite dimensional over $\mathfrak b_{\0}$. Hence the restriction of $M$
to $\mathfrak g_{\0}$ belongs to the analogue
of category $\mathcal O_n$ for $\mathfrak g_{\0}$. It is well-known
all such $\mathfrak g_{\0}$-modules have a composition series.
\end{proof}

For each $\lambda \in \Z^n$, we define the {\em Verma supermodule}
\begin{equation}
M(\lambda) := U(\mathfrak g) \otimes_{U(\mathfrak b)} \u(\lambda) \in \mathcal O_n,
\end{equation}
where $\u(\lambda)$ is as defined in the previous subsection
but viewed now as a $\mathfrak b$-supermodule by inflation from $\mathfrak h$.
By the PBW theorem, we have that
\begin{equation}\label{chv}
\operatorname{ch} M(\lambda) = 2^{\lfloor (h(\lambda)+1)/2\rfloor}
x^\lambda \prod_{1 \leq i < j \leq n}
\frac{1 + x_i^{-1} x_j}{1-x_i^{-1} x_j}.
\end{equation}
Also note by 
its definition as an induced supermodule that $M(\lambda)$ is universal
amongst all $\mathfrak g$-supermodules generated by a $\mathfrak b$-stable
submodule isomorphism to $\u(\lambda)$.
The following theorem is quite standard and parallels Theorem~\ref{class}
above.

\begin{Theorem}\label{class2}
For every $\lambda \in\Z^n$, 
$M(\la)$ has a unique irreducible quotient denoted $L(\lambda)$.
Moreover, 
\begin{itemize}
\item[(i)] $\{L(\lambda)\}_{\lambda \in \Z^n}$ is a complete set of pairwise
non-isomorphic irreducibles in category $\mathcal O_n$;
\item[(ii)]
$L(\lambda)$ is of the same type as $\u(\lambda)$,
i.e. type $\Mtype$ if $h(\lambda)$ is even, type $\Qtype$ otherwise;
\item[(iii)] 
$\ch L(\lambda) = 2^{\lfloor (h(\lambda)+1)/2 \rfloor} x^\lambda + (*)$
where $(*)$ is a linear combination of $x^\mu$ for $\mu < \lambda$;
\item[(iv)] $L(\lambda)$ is finite dimensional
if and only if $\lambda \in \Z^n_+$, in which case it is the same as
the supermodule denoted $L(\lambda)$ before.
\end{itemize}
\end{Theorem}

Let $K(\mathcal O_n)$ be the Grothendieck group of the
category $\mathcal O_n$. By Lemma~\ref{filt} and Theorem~\ref{class2}(i), 
$K(\mathcal O_n)$ 
is the free abelian group
on basis $\{[L(\lambda)]\}_{\lambda \in \Z^n}$. 
However the $[M(\la)]$'s do not form a basis for $K(\mathcal O_n)$.
To get round this problem, 
let $\widehat{K}(\mathcal O_n)$ be the completion of
$K(\mathcal O_n)$ 
with respect to the descending filtration $(K_d(\mathcal O_n))_{d \in \Z}$,
where 
$K_d(\mathcal O_n)$ is the subgroup of $K(\mathcal O_n)$ generated
by $\{[L(\lambda)]\}$ for $\la\in\Z^n$ with
$\sum_{i=1}^n i \lambda_i \geq d$.
Taking characters induces a well-defined map
$$
\operatorname{ch}:\widehat{K}(\mathcal O_n) \rightarrow 
\Z[[x_1^{\pm 1},\dots,x_n^{\pm 1}]],\qquad
[M] \mapsto \ch M.
$$
By Theorem~\ref{class2}(iii),
the characters of the irreducible supermodules in $\mathcal O_n$ are
linearly independent, hence the map $\operatorname{ch}$ is injective.
Moreover, $[M(\la)]$ equals $[L(\la)] + (*)$ where $(*)$ is a finite linear
combination of $[L(\mu)]$'s for $\mu < \la$. So working in 
$\widehat{K}(\mathcal O_n)$, we can also write
$[L(\la)]$ as $[M(\la)] + (\dagger)$ where $(\dagger)$ is a
possibly infinite linear combination of $[M(\mu)]$'s for $\mu < \la$.
This shows that the elements $\{[M(\la)]\}_{\la \in \Z^n}$
form a topological basis for $\widehat{K}(\mathcal O_n)$.

\Point{Central characters}
Let $Z$ denote the even center of the universal enveloping superalgebra
$U(\mathfrak g)$. Sergeev \cite{SergCent} has constructed an explicit set
of generators of $Z$.
For each $\lambda \in \Z^n$, let $\chi_\lambda$ be the central character
afforded by the Verma supermodule $M(\lambda)$, 
so $z \in Z$ 
acts on $M(\lambda)$ as scalar multiplication by $\chi_\lambda(z)$.
Also recall the definition of $\wt(\lambda)$
from (\ref{wtdef}), which is an element of the weight lattice
$P$ of type $\mathfrak{b}_{\infty}$ as defined in 
\ref{qg}.
As a consequence of Sergeev's description of $Z$, see e.g.
\cite[Lemma 8.9(ii),(iv)]{BKI},
we have the following fundamental fact:

\begin{Theorem}\label{hc}
For $\lambda,\mu \in \Z^n$, $\chi_\lambda = \chi_\mu$ if and only if
$\wt(\lambda) = \wt(\mu)$.
\end{Theorem}

For each $\gamma \in P$, let $\mathcal O_\gamma$
denote the full subcategory of $\mathcal O_n$
consisting of the objects all of whose irreducible subquotients
are of the form $L(\lambda)$ for $\lambda \in \Z^n$
with $\wt(\lambda) = \gamma$.
As a consequence of Theorem~\ref{hc}, we have the {\em block decomposition}
\begin{equation}\label{oblock}
\mathcal O_{n} = \bigoplus_{\gamma \in P} \mathcal O_\gamma.
\end{equation}
Moreover, $\mathcal O_\gamma$ is non-zero if and only if
$\gamma$ is a non-trivial weight of 
the tensor space $\T^{n}$ of \ref{ts}.
Note that given $\la \in \Z^n$ with $\wt(\la)=\gamma$,
 $h(\la)\equiv(\gamma,\gamma)/2 \pmod{2}$.
Hence, recalling Theorem~\ref{class2}(ii),  
all the irreducible supermodules belonging to the block $\mathcal O_\gamma$
are of the same type. We refer to this as the {\em type} of
the block $\mathcal O_\gamma$: type $\Mtype$ if $(\gamma,\gamma)/2$
is even, type $\Qtype$ if $(\gamma,\gamma)/2$ is odd.

The block decomposition (\ref{oblock}) induces the block decomposition
\begin{align}\label{kb1}
K(\mathcal O_n) &= \bigoplus_{\gamma \in P} K(\mathcal O_\gamma)
\end{align}
of the Grothendieck group, so here
$K(\mathcal O_\gamma)$ is the Grothendieck group of the
category $\mathcal O_\gamma$.
Let $\widehat{K}(\mathcal O_\gamma)$ be the closure 
of $K(\mathcal O_\gamma)$
in $\widehat{K}(\mathcal O_n)$. The elements $\{[M(\lambda)]\}$ 
for $\la \in \Z^n$
with $\wt(\lambda) = \gamma$ form a topological basis
for $\widehat{K}(\mathcal O_\gamma)$.

In a similar fashion, we have the block decomposition of the
category $\mathcal C_n$ of finite dimensional $G$-supermodules:
\begin{equation}\label{blockdec}
\mathcal C_n = \bigoplus_{\gamma \in P} \mathcal C_\gamma.
\end{equation}
This time, $\mathcal C_\gamma$ is non-zero if and only if $\gamma$ is
a non-trivial weight of the dual exterior power $\E^n$ from \ref{dcb}.
Note moreover that the natural embedding of $\mathcal C_n$
into $\mathcal O_n$ embeds $\mathcal C_\gamma$ into $\mathcal O_\gamma$.
We also get the block decomposition of the
Grothendieck group:
\begin{equation}\label{kb2}
K(\mathcal C_n) = \bigoplus_{\gamma \in P} K(\mathcal C_\gamma).
\end{equation}
The elements 
$\{[L(\lambda)]\}$ for
$\la \in \Z^n_+$ with $\wt(\la) = \gamma$
form a basis for $K(\mathcal C_\gamma)$.

\begin{Lemma}\label{trick}
$\prod_{1 \leq i < j \leq n} \frac{1 -x_i^{-1} x_j}{1+x_i^{-1} x_j}
\in 
\Z[[x_i^{-1} x_j\:|\:1 \leq i < j \leq n]]$
can be expressed
as an infinite
linear combination of $x^\mu$'s for $\mu \leq 0$ with $\wt(\mu) = 0$.
\end{Lemma}

\begin{proof}
Since the $\{[M(\mu)]\}_{\mu \in \Z^n, \wt(\mu) = 0}$ form a topological basis
for $\widehat{K}(\mathcal O_0)$, we can write
$[L(0)] = \sum_{\mu \leq 0, \wt(\mu) = 0} a_\mu [M(\mu)]$
for some coefficients $a_\mu \in \Z$.
Taking characters using (\ref{chv}) gives
$$
1 = \sum_{\substack{\mu \leq 0 \\ \wt(\mu) = 0}} 
\left(
2^{\lfloor (h(\mu)+1)/2\rfloor}a_\mu 
x^\mu
\prod_{1 \leq i < j \leq n} \frac{1+x_i^{-1}x_j}{1-x_i^{-1}x_j}\right).
$$
The lemma follows.
\end{proof}

\begin{Theorem}\label{blockthm}
For $\la \in \Z^n_+$ with $\wt(\la) = \gamma$, 
the Euler characteristic $[E(\la)]$ belongs to the block
$K(\mathcal C_\gamma)$.
\end{Theorem}

\begin{proof}
Take $\la \in \Z^n_+$. 
By Theorem~\ref{ecdef} and (\ref{pdef}), we have that
\begin{align*}
\ch E(\la) &=
2^{\lfloor(h(\la)+1)/2\rfloor}
\!\!\sum_{w \in S_n / S_\lambda}
w \left(
x^\lambda 
\prod_{\substack{1 \leq i < j \leq n\\\la_i = \la_j}} 
\frac{1 - x_i^{-1} x_j}{1 + x_i^{-1} x_j}
\prod_{1 \leq i < j \leq n}
\frac{1 + x_i^{-1} x_j}{1 - x_i^{-1} x_j}
\right)\\
&=2^{\lfloor(h(\la)+1)/2\rfloor}
\!\!\sum_{w \in S_n / S_\lambda}
\!\!(-1)^{\ell(w)} 
x^{w\lambda} 
\!\!\!\prod_{\substack{1 \leq i < j \leq n\\\la_i = \la_j}} 
\frac{1 - x_{wi}^{-1} x_{wj}}{1 + x_{wi}^{-1} x_{wj}}
\prod_{1 \leq i < j \leq n}
\frac{1 + x_i^{-1} x_j}{1 - x_i^{-1} x_j}.
\end{align*}
By Lemma~\ref{trick},
$x^{w\lambda} 
\prod_{1 \leq i < j \leq n, \la_i = \la_j} 
\frac{1 - x_{wi}^{-1} x_{wj}}{1 + x_{wi}^{-1} x_{wj}}$
is a (possibly infinite) linear combination 
of $x^\mu$'s for $\mu \leq w\lambda$ with $\wt(\mu) = 
\wt(w \lambda) = \gamma$.
Hence recalling (\ref{chv}), $\ch E(\la)$ is a 
(possibly infinite) linear
combination of $\ch M(\mu)$'s for 
$\mu\in\Z^n$ with $\wt(\mu) = \gamma$.
Working instead in 
$\widehat{K}(\mathcal O_n)$, 
we have shown that we can write
$[E(\la)]$ as a (possibly infinite)
linear combination of $[M(\mu)]$'s for $\mu \in \Z^n$
with $\wt(\mu) = \gamma$, i.e.
$[E(\la)] \in \widehat{K}(\mathcal O_\gamma)$.
The theorem follows.
\end{proof}

The theorem immediately implies that the
decomposition numbers $d_{\mu,\la}$ defined in (\ref{split})
and the inverse decomposition numbers $d_{\mu,\la}^{-1}$
defined in (\ref{split2}) are zero 
whenever $\wt(\mu) \neq \wt(\la)$. In particular:

\begin{Corollary}\label{cory}
The elements $\{[E(\la)]\}$ for
$\la \in \Z^n_+$ with $\wt(\la) = \gamma$ form a basis
for $K(\mathcal C_\gamma)$. 
\end{Corollary}

\Point{Translation functors}\label{tf}
Now we are ready to link the Grothendieck group $K(\mathcal C_n)$ 
with the $\U_\Z$-module $\E^n_{\Z}$ constructed
in \ref{zforms}.
We define an isomorphism of abelian groups
\begin{equation}\label{iotadef}
\iota:K(\mathcal C_n) \rightarrow \E^n_{\Z},\qquad
[E(\la)] \mapsto E_\la(1)\qquad(\la \in \Z_n^+).
\end{equation}
Using $\iota$, we lift the actions of the
generators of $\U_{\Z}$
on $\E^n_{\Z}$ from (\ref{generators}) to define an action
of $\U_{\Z}$ directly on the Grothendieck group $K(\mathcal C_n)$.
Note Corollary~\ref{cory} shows that the block decomposition (\ref{kb2})
coincides with the usual weight space decomposition 
of $K(\mathcal C_n)$ as a $\U_\Z$-module.
Using Lemma~\ref{expact} specialized at $q = 1$, we can write down
the action of the operators
$E_i, F_i$ on $K(\mathcal C_n)$ explicitly:
\begin{align}\label{formulas1}
E_i [E(\lambda)] &= 
\sum_{\substack{r\,\text{with}\,\lambda-\delta_r\in\Z^n_+,\\\sigma_r = -,-+\,\text{or}\,--}}
c_{\la,r} [E(\lambda-\delta_r)],\\
F_i [E(\lambda)] &=\label{formulas2}
\sum_{\substack{r\,\text{with}\,\lambda+\delta_r\in\Z^n_+,\\\sigma_r = +,-+\,\text{or}\,++}}
c_{\la,r} [E(\lambda+\delta_r)],
\end{align}
for $\la \in \Z^n_+$ and $i \geq 0$. 
Here,  $(\sigma_1,\dots,\sigma_n)$ is the $i$-signature
of $\lambda$ defined according to (\ref{isig}) and
$$
c_{\la,r} = \left\{
\begin{array}{ll}
1+(-1)^{z(\lambda)}
&\hbox{if $\sigma_r = --$ or $++$,}\\
1&\hbox{otherwise.}
\end{array}\right.
$$
The goal in the remainder of this subsection is to give a representation
theoretic interpretation of the operators $E_i$ and $F_i$.

We need certain ``translation functors''
$$
\Tr^i, \Tr_i:\mathcal O_n\rightarrow \mathcal O_n,
$$
one for each $i \geq 0$, similar to the functors
defined by Penkov and Serganova in \cite[$\S$2.3]{PS1}.
Recalling the decomposition (\ref{blockdec}),
it suffices to define the functors $\Tr^i, \Tr_i$ on the subcategory 
$\mathcal O_\gamma$
for some fixed $\gamma \in P$, since one can then extend
additively to get the definition on $\mathcal O_n$ itself.
We will write
$\pr_\gamma:\mathcal O_n \rightarrow \mathcal O_\gamma$
for the natural projection functor.
Given $M \in \mathcal O_\gamma$, let
\begin{equation}
\Tr^i M = \pr_{\gamma - \alpha_i}(M \otimes V),
\qquad
\Tr_i M = \pr_{\gamma + \alpha_i}(M \otimes V^*).
\end{equation}
On a morphism $f:M \rightarrow N$ in $\mathcal O_\gamma$,
$\Tr^i f$ and $\Tr_i f$ are the restrictions of the
maps $f \otimes \operatorname{id}$.
Obviously these functors send finite dimensional
supermodules to finite dimensional supermodules,
so also give us functors
$$
\Tr^i, \Tr_i:\mathcal C_n \rightarrow \mathcal C_n
$$
by restriction.

\begin{Lemma}\label{adhjs} For each $i \geq 0$, the
functors $\Tr^i$ and $\Tr_i$ (on either of the categories
$\mathcal O_n$ or $\mathcal C_n$)
are both left and right adjoint to each
other, hence both are exact functors.
Moreover, there are natural isomorphisms
\begin{align}
(\Tr^i M)^* &\cong \Tr_i (M^*),\label{stardual}
&(\Tr_i M)^* &\cong \Tr^i (M^*),
\\
(\Tr^i M)^\tau &\cong \Tr^i (M^\tau)
&
(\Tr_i M)^\tau &\cong \Tr_i (M^\tau)\label{taudual}
\end{align}
for each $M \in \mathcal C_n$.
\end{Lemma}

\begin{proof}
The first part is a well-known fact about translation functors.
For the second part, we obviously have natural isomorphisms
$(M \otimes V)^* \cong M^* \otimes V^*$,
$(M \otimes V)^\tau \cong M^\tau \otimes V$.
So it suffices to show that $*$-duality maps 
$\mathcal C_\gamma$ into $\mathcal C_{-\gamma}$.
and that $\tau$-duality maps $\mathcal C_\gamma$ into itself.
This is immediate from Theorem~\ref{class}(iv).
\end{proof}

\begin{Lemma}\label{vmod} For $\la \in \Z^n_+$ and $i \geq 0$, and
a finite dimensional $G$-supermodule $M$ belonging to the block
$\mathcal C_\gamma$ for some $\gamma \in P$,
we have that
\begin{align*}
[\Tr_i M] = \left\{\begin{array}{ll}
E_i [M] &\hbox{if $i = 0$ and $\mathcal C_\gamma$ is of type $\Mtype$,}\\
2 E_i [M] & \hbox{if $i \neq 0$ or $\mathcal C_\gamma$ is of type $\Qtype$;}
\end{array}\right.\\
[\Tr^i M] = \left\{\begin{array}{ll}
F_i [M] &\hbox{if $i = 0$ and $\mathcal C_\gamma$ is of type $\Mtype$,}\\
2 F_i [M] & \hbox{if $i \neq 0$ or $\mathcal C\gamma$ is of type $\Qtype$.}
\end{array}\right.
\end{align*}
\end{Lemma}

\begin{proof}
We explain the argument for $\Tr^i$, the case of $\Tr_i$ being similar.
Since the functor $\Tr^i$ is exact, it induces an additive
map also denoted $\Tr^i$ on $K(\mathcal C_n)$.
By Corollary~\ref{cory}, the $\{[E(\la)]\}$ for $\la \in \Z^n_+$
with $\wt(\la) = \gamma$ form a basis for $K(\mathcal C_\gamma)$.
Therefore it suffices to consider
$\Tr^i [E(\la)]$ for such a $\la$.

Observe that $\ch V = 2 (x_1+\dots + x_n)$ and that 
the functor $?\otimes V$ is isomorphic to the functor
$\bigoplus_{i \geq 0} \Tr^i$.
So a calculation using Theorem~\ref{ecdef}
and Lemma~\ref{prodform} gives that
$$
\bigoplus_{i \geq 0} \Tr^i [E_\la]
=
\sum_{r \,\text{with}\, \la+\delta_r\in\Z^n_+}
d_{\la,r} [E_{\la+\delta_r}]
$$
where
$$
d_{\la,r} = \left\{\begin{array}{ll}
2&\hbox{if $\la_r \neq -1,0$,}\\
1&\hbox{if $\la_r = 0$ and $h(\la)$ is even,}\\
2&\hbox{if $\la_r = 0$ and $h(\la)$ is odd,}\\
0&\hbox{if $\la_r = -1$ and $z(\la)$ is odd,}\\
2&\hbox{if $\la_r = -1$, $z(\la)$ is even and $h(\la)$ is even,}\\
4&\hbox{if $\la_r = -1$, $z(\la)$ is even and $h(\la)$ is odd.}
\end{array}\right.
$$
Now apply $\pr_{\wt(\la)-\alpha_i}$ to both sides and use
Theorem~\ref{blockthm} and (\ref{formulas2}) to deduce that
$$
\Tr^i [E(\la)] = 
\left\{\begin{array}{ll}
F_i [E(\la)] &\hbox{if $i = 0$ and $h(\la)$ is even,}\\
2 F_i [E(\la)] & 
\hbox{if $i \neq 0$ or $h(\la)$ is odd.}
\end{array}\right.
$$
This completes the proof.
\end{proof}

\Point{\boldmath Crystal structure}
In this section we relate the structure 
of the supermodules $\Tr^i L(\la)$ and $\Tr_i L(\la)$ for 
certain weights $\la$
to the crystal operators defined in \ref{cs}.
The methods employed here are essentially the same as those
of Penkov and Serganova in \cite[$\S$2.3]{PS1}.
Throughout the subsection we will use the following notation:
for supermodules $X$ and $Y$, $\frac{X}{Y}$ will denote some
extension of $X$ by $Y$,
and $[M:L]$ will denote the composition multiplicity of an irreducible
$L$ in $M$.

\begin{Lemma}\label{calc}
Let $\la = (1,0^{n-2}, -1)$.
Then $[M(\la):L(0)] \geq 2$.
\end{Lemma}

\begin{proof}
Pick a basis $v,v'$ for $M(\la)_\la$ such that
$h_1' v = v', h_1' v' = v, h_n' v = v', h_n' v' = -v$ and
$h_r' v= h_r' v' = 0$ for $r = 2,\dots,n-1$.
Let
$$
x = e_{2,1}' (e_{3,2}' + e_{3,2} h_3') \dots
(e_{n,n-1}' + e_{n,n-1} h_n') v, \qquad y = h_1' x.
$$
We claim that $\langle x, y \rangle$ is a two dimensional 
indecomposable $\mathfrak b$-submodule of $M(\la)$.
The lemma follows immediately form this, 
since both $x$ and $y$ are of weight $0$.
To prove the claim, we proceed by induction on $n$. The base case $n=2$
follows on checking that $x = e_{2,1}' v$ is annihilated by
$e_{1,2}$ and $e_{1,2}'$ and that $y = h_1' x = h_2' x$.
For $n > 2$, we set
\begin{align*}
w&=(e_{n,n-1}'+e_{n,n-1} h_n') v = e_{n,n-1}' v + e_{n,n-1} v',\\
w'&=-(e_{n,n-1}'+e_{n,n-1} h_n') v' = e_{n,n-1} v - e_{n,n-1}'v'.
\end{align*}
Now one checks that
$h_1' w = w', h_1' w' = w, h_{n-1}' w = w', h_{n-1}' w' = -w$ and
that $h_r' w = h_r' w' = 0$ for each $r =2,\dots,n-2,n$.
Moreover, $w,w'$ are annihilated by 
$e_{r,r+1},e_{r,r+1}'$ for $r = 1,\dots,n-1$.
Now the result follows from the induction hypothesis.
\end{proof}

\begin{Lemma}\label{ps}
Let $\la \in \Z_+^n$ with $\la_r > 0$ and $\la_r+\la_s = 0$ 
for some $1 \leq r < s \leq n$.
Then,
$$
[M(\la):L(\la - \delta_r+\delta_s)]
=\left\{
\begin{array}{ll}
1&\hbox{if $\la_r > 1$,}\\
2&\hbox{if $\la_r = 1$.}
\end{array}\right.
$$
\end{Lemma}

\begin{proof}
We will actually only prove here that 
the left hand side is $\geq$ the right hand side, since that is all
we really need later on. The fact that the left hand side actually {\em equals}
the right hand side can easily 
be extracted from the proof of Lemma~\ref{l2} below.
First suppose that $\la_r >1$.
By \cite[Proposition 2.1]{PS1},
$$
\hom_{\mathfrak g} (M(\la-\delta_r+\delta_s), M(\la)) \neq 0.
$$
This immediately implies that $[M(\la):L(\la-\delta_r+\delta_s)] \geq 1$.
Now assume that $\la_r = 1$, when we need to show that in fact
$[M(\la):L(\la-\delta_r+\delta_s)] \geq 2$.
Let $\mathfrak{p}$ be the upper triangular parabolic subalgebra
of $\mathfrak g$
of type $(1^{r-1}, s-r+1, 1^{n-s})$. Consider the
Verma supermodule $M_{\mathfrak p}(\la) = U(\mathfrak p)
\otimes_{U(\mathfrak b)} \u(\la)$
over $\mathfrak p$.
By Lemma~\ref{calc}, it contains $L_{\mathfrak p}(\la-\delta_r+\delta_s)$
with multiplicity at least two. 
Since
$U(\mathfrak g) \otimes_{U(\mathfrak p)} M_{\mathfrak p}(\la) \cong M(\la)$
and 
$U(\mathfrak g) \otimes_{U(\mathfrak p)} 
L_{\mathfrak p}(\la-\delta_r+\delta_s) \twoheadrightarrow L(\la-\delta_r+\delta_s)$, the lemma follows.
\end{proof}

\begin{Lemma}
Let $\la \in \Z^n$. Then $M := M(\la) \otimes V$ has a filtration
$0 = M_0 < M_1 < \dots < M_n = M$ such that
$$
M_r / M_{r-1} \cong
\left\{
\begin{array}{ll}
M(\la+\delta_r) \oplus\Pi M(\la+\delta_r)&\hbox{if $\la_r \neq -1,0$}\\
M(\la+\delta_r)&\hbox{if $\la_r = 0$ and $h(\la)$  even,}\\
M(\la+\delta_r)\oplus\Pi M(\la+\delta_r)&\hbox{if $\la_r = 0$ and $h(\la)$  odd,}\\
\frac{M(\la+\delta_r)}{M(\la+\delta_r)}&\hbox{if $\la_r = -1$ and $h(\la)$ even,}\\
\frac{\Pi M(\la+\delta_r)}{M(\la+\delta_r)}
\oplus \frac{M(\la+\delta_r)}{\Pi M(\la+\delta_r)}
 &\hbox{if $\la_r = -1$ and $h(\la)$ odd,}
\end{array}
\right.
$$
for each $r = 1,\dots,n$.
\end{Lemma}

\begin{proof}
Note $V \cong \u(\delta_1) \oplus\dots\oplus \u(\delta_n)$
as an $\mathfrak h$-supermodule.
Now the lemma follows by a standard construction, see e.g.
\cite[Lemma 4.24]{KL} for a similar situation, from the observation
that
$$
\u(\la) \otimes \u(\delta_r)
\cong
\left\{
\begin{array}{ll}
\u(\la+\delta_r) \oplus\Pi \u(\la+\delta_r)&\hbox{if $\la_r \neq -1,0$}\\
\u(\la+\delta_r)&\hbox{if $\la_r = 0$ and $h(\la)$  even,}\\
\u(\la+\delta_r)\oplus\Pi \u(\la+\delta_r)&\hbox{if $\la_r = 0$ and $h(\la)$  odd,}\\
\frac{\u(\la+\delta_r)}{\u(\la+\delta_r)}&\hbox{if $\la_r = -1$ and $h(\la)$ even,}\\
\frac{\Pi \u(\la+\delta_r)}{\u(\la+\delta_r)}
\oplus \frac{\u(\la+\delta_r)}{\Pi \u(\la+\delta_r)}
 &\hbox{if $\la_r = -1$ and $h(\la)$ odd}
\end{array}
\right.
$$
as $\mathfrak h$-supermodules.
We prove this in just two of the situations, the rest being
similar.

First suppose $\la_r = -1$ and $h(\la)$ is odd.
By character considerations $\u(\la) \otimes \u(\delta_r)$
has just four composition factors, all isomorphic to $\u(\la+\delta_r)$.
Moreover, both $\u(\la)$ and $\u(\delta_r)$ possess odd automorphisms
$J_1, J_2$ respectively with $J_i^2 = -1$. 
The $\pm \sqrt{-1}$-eigenspaces of the even automorphism
$J_1 \otimes J_2$ of $\u(\la) \otimes \u(\delta_r)$ decompose it into
a direct sum of two factors, and the map $J_1 \otimes 1$ gives an odd
isomorphism between the factors.
Hence $\u(\la)\otimes\u(\delta_r)$ has the structure given.

Second suppose $\la_r \neq -1,0$ and $h(\la)$ is even.
In that case, $\u(\la)$ is of type $\Mtype$, and by character
considerations $\u(\la) \otimes \u(\delta_r)$ has just two composition
factors, both isomorphic to $\u(\la+\delta_r)$ which also has type $\Mtype$.
Let $J$ be an odd automorphism of $\u(\delta_r)$.
Then, $1 \otimes J$ is an odd automorphism of $\u(\la) \otimes \u(\delta_r)$.
If $\u(\la) \otimes \u(\delta_r)$ was to be a non-split extension of
$\u(\la+\delta_r)$ and $\Pi \u(\la+\delta_r)$, it would possess a nilpotent
odd endomorphism, but no odd automorphism. Hence it must split as a
direct sum as required.
\end{proof}

Note on applying the exact functor $\pr_{\wt(\la)-\alpha_i}$ to the filtration
constructed in the lemma, we obtain a filtration of
$\Tr^i M(\la)$. Moreover, if $\la \in \Z^n_+$, $L(\la)$ is a finite dimensional
quotient of $M(\la)$, so $\Tr^i L(\la)$ is a finite dimensional quotient
of $\Tr^i M(\la)$. Hence we also get induced a filtration of
$\Tr^i L(\la)$ 
whose the factors are finite dimensional quotients of the factors
of the filtration of $\Tr^i M(\la)$.
We refer to this as the {\em canonical filtration} of $\Tr^i L(\la)$.
Its properties are used repeatedly in the proofs of the next two lemmas.

\begin{Lemma}\label{l1} Suppose $\la \in \Z^n_+$ and $i \geq 0$ are 
such that
$\phi_i(\la) = 1$. 
Then, $F_i [L(\la)] = [L(\tilde F_i (\la))]$.
\end{Lemma}

\begin{proof}
We actually prove the slightly stronger statement that
$$
\Tr^i L(\la) \cong
\left\{
\begin{array}{ll}
L(\mu)&\hbox{if $i = 0$ and $h(\la)$ is even,}\\
L(\mu) \oplus \Pi L(\mu)&\hbox{if $i > 0$ or $h(\la)$ is odd,}
\end{array}
\right.
$$
where $\mu := \tilde F_i(\la)$.
The lemma follows from this and Lemma~\ref{vmod} above.
The possible cases for $\la, \mu$ are 
listed explicitly in \ref{cs}, in particular we have that
$\mu = \la + \delta_s$ for some $1 \leq s \leq n$.
By considering the canonical filtration in each of the cases,
$\Tr^i M(\la) \cong M(\mu)$
if $i = 0$ and $h(\la)$ is even,
$\Tr^i M(\la) \cong M(\mu) \oplus \Pi M(\mu)$ otherwise.
In particular, $[\Tr^i M(\la):L(\mu)] = 1$
if $i =  0$ and $h(\la)$ is even, $[\Tr^i M(\la):L(\mu)] = 2$ otherwise.
Now, $\Tr^i L(\la)$ is a quotient of $\Tr^i M(\la)$, and moreover it
is self-dual under the duality $\tau$
by Theorem~\ref{class}(iv) and (\ref{taudual}). 
Hence, $\Tr^i L(\la)$ is necessarily a quotient of
$L(\mu)$ if $i = 0$ and $h(\la)$ is even,
$L(\mu) \oplus \Pi(\mu)$ otherwise.
So to complete the proof, we just need to show that
$[\Tr^i L(\la):L(\mu)] = [\Tr^i M(\la):L(\mu)]$.

Suppose for a contradiction that $[\Tr^i L(\la):L(\mu)] <
[\Tr^i M(\la):L(\mu)]$. Then 
there must be some composition factor $L(\la')\not\cong L(\la)$ of $M(\la)$
such that $[\Tr^i L(\la') : L(\mu)] > 0$.
Hence, $[\Tr^i M(\la') : L(\mu)] > 0$ for some $\la' < \la$
with $\wt(\la') = \wt(\la)$.
Considering the canonical filtration of
$\Tr^i M(\la')$, there must exist some $1 \leq r \leq n$ such that
$\mu = \la + \delta_s \leq \la' + \delta_r$.
But then $\la  - \delta_r + \delta_s \leq \la' < \la$
and $\wt(\la') = \wt(\la)$. 
It is easy to see in each case that there is 
no such $\la'$, giving the 
desired contradiction.
\end{proof}

\begin{Lemma}\label{l2}
Suppose $\la \in \Z^n_+$ and $i \geq 0$ are such that $\phi_i(\la) = 0$.
Then, $F_i [L(\la)]  = 0$.
\end{Lemma}

\begin{proof}
The possibilities for $\la$ are listed explictly
in \ref{cs}. In almost all of the configurations, 
we get that $\Tr^i L(\la) = 0$, hence $F_i [L(\la)]= 0$,
immediately
by looking at the canonical filtration.
There are just two difficult cases in which we need to argue further.

In the first case, $i = 0$ and $\la = (\cdots, 1, 0^{r-1}, -1, \cdots)$ for some
$r \geq 1$, where $\cdots$ denote entries different from $-1,0,1$.
We let
$\mu = (\cdots,0^{r+1},\cdots)$ and $\nu = (\cdots, 1, 0^r, \cdots)$,
where the $\cdots$ are the same entries as in $\la$.
Also let $c = 1$ if $h(\la)$ is even, $c = 2$ if $h(\la)$ is odd.
By considering the canonical filtration, we see here that
$[\Tr^i M(\la) : L(\nu)]  = 2c$.
By Lemma~\ref{ps}, $[M(\la): L(\mu)] \geq 2$,
hence $[\Tr^i M(\la): L(\nu)] \geq [\Tr^i L(\la):L(\nu)] 
+ 2 [\Tr^i L(\mu):L(\nu)]$.
By Lemma~\ref{l1}, $[\Tr^i L(\mu):L(\nu)] = c$.
This shows that $[\Tr^i L(\la):L(\nu)] = 0$ (and also 
that $[M(\la):L(\mu)] = 2$, completing the
proof of Lemma~\ref{ps} in this case).
But by the canonical filtration, $\Tr^i L(\la)$ has a filtration where
all the factors are quotients of $M(\nu)$, so it has to be zero.

In the second case,  $i > 0$ and 
$\la = (\dots, i+1, \dots, -i-1,\dots)$, 
where $\dots$ denote entries different from $-i-1,-i,i,i+1$.
We let
$\mu = (\cdots,i,\cdots,-i,\cdots)$ and 
$\nu = (\cdots, i+1, \cdots,-i,\cdots)$,
where the $\cdots$ are the same entries as in $\la$.
By the canonical filtration, 
$\Tr^i M(\la) \cong M(\nu) \oplus \Pi M(\nu)$, hence
$[\Tr^i M(\la): L(\nu)]  = 2$.
By Lemma~\ref{ps}, $[M(\la): L(\mu)] \geq 1$,
hence $[\Tr^i M(\la): L(\nu)] \geq [\Tr^i L(\la):L(\nu)] 
+ [\Tr^i L(\mu):L(\nu)]$.
By Lemma~\ref{l1}, $[\Tr^i L(\mu):L(\nu)] = 2$.
This shows that $[\Tr^i L(\la):L(\nu)] = 0$ (and also 
completes the proof of Lemma~\ref{ps} in this case).
Hence, since $\Tr^i L(\la)$ is a quotient of
$M(\nu) \oplus\Pi M(\nu)$,  it has to be zero.
\end{proof}

\Point{Injective supermodules}
We refer at this point 
to \cite[I.3]{Jantzen} for the general facts about injective
modules over a group scheme, all of which generalize to supergroups.
In particular, for every $G$-supermodule $M$,
there is an injective $G$-supermodule $U$, unique up to isomorphism,
such that $\soc_G M \cong \soc_G U$. 
We call $U$ the {\em injective hull} of $M$. 
For $\la \in \Z^n_+$, let $U(\lambda)$ denote the injective hull of
$L(\lambda)$.
Any injective $G$-supermodule $M$ is isomorphic to a direct sum
of $U(\la)$'s, the number of summands isomorphic to $U(\lambda)$
being equal to the multiplicity of $L(\lambda)$ in $\soc_G M$.

\begin{Lemma}\label{fd} Each $U(\la)$ is finite dimensional.
\end{Lemma}

\begin{proof}
Consider the functor
$\ind_{G_{\ev}}^G := \hom_{U(\mathfrak g_{\0})}(U(\mathfrak g), ?).$
Since $U(\mathfrak g)$ is a free right $U(\mathfrak g_{\0})$-supermodule
of finite rank, this is an exact functor mapping finite dimensional
$G_{\ev}$-supermodules to finite dimensional $G$-supermodules.
It is right adjoint to the restriction functor $\res^G_{G_{\ev}}$
from the category of $G$-supermodules to the category of 
$G_{\ev}$-supermodules, so sends injectives to injectives.
Now take  $\la \in \Z^n_+$.
Since {\em every} $G_{\ev}$-module
is injective, we get that
$\ind_{G_{\ev}}^G \res^G_{G_{\ev}} L(\la)$ 
is a finite dimensional injective $G$-supermodule.
Moreover, the unit of the adjunction gives an embedding of $L(\la)$ into
$\ind_{G_{\ev}}^G \res^G_{G_{\ev}} L(\la)$.
Hence the injective hull of $L(\la)$ is finite dimensional.
\end{proof}

Let $\mathcal C_n^*$ be the category of all finite dimensional injective
$G$-supermodules. 
The block decomposition (\ref{blockdec}) of $\mathcal C_n$ induces 
an analogous block decomposition of the subcategory $\mathcal C_n^*$
\begin{equation}
\mathcal C_n^* = \bigoplus_{\gamma \in P} \mathcal C_\gamma^*.
\end{equation}
Let $K(\mathcal C_n^*)$ (resp. $K(\mathcal C_\gamma^*)$)
be the Grothendieck group of the category
$\mathcal C_n^*$ (resp. $\mathcal C_\gamma^*$).
By Lemma~\ref{fd}, $K(\mathcal C_n^*)$
is the free abelian group on basis $\{[U(\lambda)]\}_{\la \in \Z^n_+}$,
and $K(\mathcal C_\gamma^*)$ is the subgroup generated
by the $\{[U(\la)]\}$ for $\la \in \Z^n_+$ with $\wt(\la) = \gamma$.

Form the completion 
$\widehat{K}(\mathcal C_n^*)$ of the Grothendieck
group $K(\mathcal C_n^*)$ with respect to the descending filtration
$(K_d(\mathcal C_n^*))_{d \in \Z}$ where
$K_d(\mathcal C_n^*)$ is the subgroup generated by
$\{[U(\lambda)]\}$ for $\lambda \in \Z^n_+$ with
$\sum_{i=1}^n i \lambda_{n+1-i} \geq d$.
The important thing in this definition is that 
vectors of the form $[U(\lambda)] + 
(*)$  make sense whenever $(*)$ is an infinite linear combination of
$[U(\mu)]$'s for $\mu > \la$.
In particular, the following are well-defined elements of
$\widehat{K}(\mathcal C_n^*)$ for each $\la \in \Z^n_+$:
\begin{equation}\label{ffdef}
[F(\lambda)] := \sum_{\mu \in \Z^n_+} d_{-w_0\mu,-w_0\la}^{-1} 
[U(\mu)],
\end{equation}
recall (\ref{split2}).
By the unitriangularity of the inverse decomposition numbers,
the elements
$\{[F(\lambda)]\}_{\la \in \Z^n_+}$ give a topological basis for 
$\widehat{K}(\mathcal C_n^*)$. Note $[F(\la)]$ does not 
in general belong to $K(\mathcal C_n^*)$ itself.

We define a pairing $\langle.,.\rangle$ between the Grothendieck
groups $K(\mathcal C_n)$ and $K(\mathcal C_n^*)$ 
by letting
\begin{equation}\label{theform}
\langle [L(-w_0\la)], [U(\mu)] \rangle = \delta_{\la,\mu}
\end{equation}
for each $\la,\mu \in \Z^n_+$. 
The pairing $\langle.,.\rangle$ extends by continuity to 
give a pairing also denoted 
$\langle.,.\rangle$ 
between $K(\mathcal C_n)$ and $\widehat{K}(\mathcal C_n^*)$.
In that case, by the definitions
(\ref{split}), (\ref{ffdef}) and (\ref{theform}), we have that
\begin{equation}\label{theform2}
\langle [E(-w_0\la)], [F(\mu)] \rangle = \delta_{\la,\mu}
\end{equation}
for each $\la,\mu \in \Z^n_+$. 
We record the following lemma which follows from 
a standard property of injective hulls, see
\cite[I.3.17(3)]{Jantzen}.

\begin{Lemma}\label{foror}
Suppose we are given $M \in \mathcal C_n$ and $U \in \mathcal C_\gamma^*$
for some $\gamma \in P$.
Then, 
\begin{equation*}
\langle [M], [U] \rangle =
\left\{
\begin{array}{ll}
\dim \hom_G(M^*, U)&\hbox{if $\mathcal C_\gamma$ is of type $\Mtype$,}\\
\frac{1}{2}\dim \hom_G(M^*, U)&\hbox{if $\mathcal C_\gamma$ is of type $\Qtype$.}
\end{array}
\right.
\end{equation*}
\end{Lemma}

Now recall the definition of the $\U_\Z$-module $\FF^n_{\Z}$ from
\ref{zforms}. 
We define a continuous isomorphism
\begin{equation}\label{iotastardef}
\iota^*:\FF^n_{\Z} \rightarrow \widehat{K}(\mathcal C^*_n),\qquad
F_\la(1) \mapsto [F(\la)]\qquad(\la \in \Z_n^+).
\end{equation}
The following lemma shows that this map  
$\iota^*$ is the dual map to $\iota$ from (\ref{iotadef}) with respect to 
the pairings $\langle.,.\rangle$:

\begin{Lemma}\label{adj}
$\langle \iota([M]), v \rangle = \langle [M], \iota^* (v) \rangle$
for all $[M] \in K(\mathcal C_n)$ and $v \in \FF^n_{\Z}$.
\end{Lemma}

\begin{proof}
It suffices to check this for $[M] = [E(-w_0 \lambda)]$ and $v = F_\mu(1)$
for $\la,\mu \in \Z^n_+$, when it follows immediately from
(\ref{formy}) and (\ref{theform2}).
\end{proof}

As in \ref{tf}, we lift the action of $\U_\Z$
on 
$\FF^n_{\Z}$ to the completed Grothendieck group
$\widehat{K}(\mathcal C_n^*)$ through the isomorphism $\iota^*$.
By Lemmas~\ref{adj} and \ref{sym}, we have at once that
\begin{equation}\label{adjs}
\langle E_i [L], [U] \rangle = \langle [L], E_i [U] \rangle,
\qquad
\langle F_i [L], [U] \rangle = \langle [L], F_i [U] \rangle,
\end{equation}
for each $i \geq 0$ and all
$[L] \in K(\mathcal C_n), [U] \in \widehat{K}(\mathcal C_n^*)$.
The next lemma gives a purely representation theoretic
interpretation of the operators $E_i$ and $F_i$
on $\widehat{K}(\mathcal C_n^*)$, cf. Lemma~\ref{vmod}.
For the statement, recall by Lemma~\ref{adhjs} that
the functors $\Tr_i$ and $\Tr^i$ 
send injectives to injectives.

\begin{Lemma}
For $\la \in \Z^n_+$ and $i \geq 0$, there exist 
injective $G$-supermodules $E_i U(\la)$ and $F_iU(\la)$,
unique up to isomorphism, characterized by
\begin{align*}
\Tr_i U(\la) &\cong \left\{\begin{array}{ll}
E_i U(\la) &\hbox{if $i = 0$ and $h(\la)$ is even,}\\
E_i U(\la) \oplus \Pi E_i U(\la) & \hbox{if $i \neq 0$ or $h(\la)$ is odd;}
\end{array}\right.\\
\Tr^i U(\la) &\cong \left\{\begin{array}{ll}
F_i U(\la) &\hbox{if $i = 0$ and $h(\la)$ is even,}\\
F_i U(\la)\oplus \Pi F_i U(\la) &\hbox{if $i \neq 0$ or $h(\la)$ is odd.}
\end{array}\right.
\end{align*}
Moreover, $E_i [U(\la)]  = [E_i U(\la)]$ and $F_i [U(\la)]  = [F_i U(\la)]$.
\end{Lemma}

\begin{proof}
We just explain the proof for $F_i$.
Take $\la \in \Z^n_+$ and $i \geq 0$.
Note uniqueness of $F_i U(\la)$ is immediate by Krull-Schmidt.
For existence, we consider three separate cases.

{\em Case one.} If $i = 0$ and $h(\la)$ is even, let
$F_i U(\la) := \Tr^i U(\la)$.

{\em Case two.} If $h(\la)$ is odd,
$L(\la)$ is of type $\Qtype$, so possesses an odd automorphism
$J_1$ with $J_1^2 = -1$. 
This induces an odd automorphism also denoted $J_1$ of the
injective hull $U(\la)$. Also the natural representation $V$ possesses
the odd automorphism $J$ defined earlier (\ref{jdef}).
The map $J_1 \otimes J$ induces an even automorphism
of the summand $\Tr^i U(\la)$ of $U(\la) \otimes V$.
Its $\pm \sqrt{-1}$-eigenspaces decompose
$\Tr^i U(\la)$ into a direct sum of
two $G$-supermodules, and the map $1 \otimes J$ is an odd isomorphism
between them. Let $F_i U(\la)$ be the $\sqrt{-1}$-eigenspace (say),
then $\Tr^i U(\la) \cong F_i U(\la) \oplus \Pi F_i U(\la)$ as required.

{\em Case three.} If $h(\la)$ is even and $i > 0$, then
the map $1 \otimes J$ induces an odd automorphism of $\Tr^i U(\la)$,
hence also of the socle $S$ of $\Tr^i U(\la)$.
All constituents of $S$ are of type $\Mtype$, so $S$ must decompose
as a $G$-supermodule as $S_- \oplus S_+$ with $(1 \otimes J) S_{\pm} = 
S_{\mp}$.
This decomposition of the socle induces a decomposition of
the injective supermodule
$\Tr^i U(\la)$, say $\Tr^i U(\la) = M_- \oplus M_+$
where $\soc_G M_{\pm} = S_{\pm}$.
In this case we let $F_i U(\la) = M_+$ (say).

It remains to show that $F_i [U(\la)]  = [F_i U(\la)]$. 
To do this, it suffices to prove that
$\langle [L(\mu)], F_i[U(\la)] \rangle = 
\langle [L(\mu)], [F_i U(\la)] \rangle$
for all $\mu \in \Z^n_+$.
This is done using Lemmas~\ref{foror}, \ref{vmod},
(\ref{adjs}) and the adjointness
of $\Tr^i, \Tr_i$ from Lemma~\ref{adhjs}.
We just explain the argument in the case that $h(\la)$ is even and $i = 0$,
the other situations being entirely similar.
First note by weight considerations
that both sides of the identity we are trying to verify
are zero unless $h(\mu)$ is odd. Now compute:
\begin{align*}
\langle [L(\mu)], F_i [U(\la)] \rangle &=
\langle F_i [L(\mu)],  [U(\la)] \rangle
= \frac{1}{2} \langle [\Tr^i L(\mu)], [U(\la)] \rangle\\
&=
\frac{1}{2} \dim \hom_G((\Tr^i L(\mu))^*, U(\la))\\&
=
\frac{1}{2} \dim \hom_G(\Tr_i (L(\mu)^*), U(\la))\\
&=
\frac{1}{2} \dim \hom_G(L(\mu)^*, \Tr^i U(\la))\\
&=
\langle L(\mu), [\Tr^i U(\la)] \rangle
=\langle L(\mu), [F_i U(\la)] \rangle.
\end{align*}
This completes the proof.
\end{proof}

\begin{Lemma}\label{cor2} Let $\la \in \Z^{n}_+$ and $i \geq 0$.
\begin{itemize}
\item[(i)]
If $\phi_i(\la) = 0$ and $\eps_i(\la) > 0$ then
$E_i U(\lambda) \cong
U(\tilde E_i(\la))$.
\item[(ii)]
If $\eps_i(\la) = 0$ and $\phi_i(\la) > 0$ then
$F_i U(\lambda) = 
U(\tilde F_i (\la))$.
\end{itemize}
\end{Lemma}

\begin{proof}
We just prove (ii), since (i) follows by applying $*$.
We know that $F_i U(\la)$ is a direct sum of injective indecomposables.
To compute the multiplicity of $U(\mu)$ for a given $\mu \in \Z^n_+$,
it suffices to compute
$$
\langle [L(\mu)], [F_i U(\la)] \rangle =
\langle [L(\mu)], F_i [U(\la)] \rangle =
\langle F_i [L(\mu)], [U(\la)] \rangle.
$$
By block considerations and Lemmas~\ref{l1} and \ref{l2}, 
that is zero unless $\mu = -w_0 \tilde F_i(\lambda)$,
in which case it is one. Hence
$F_i U(\la) = U(\tilde F_i(\lambda))$.
\end{proof}

Now we can construct the injective supermodules $U(\la)$.
If $\la$ is typical, there are no problems:

\begin{Lemma}\label{typ} Suppose that $\la \in \Z^n_+$ is typical.
Then $U(\la) = H^0(\la) = L(\la)$.
\end{Lemma}

\begin{proof}
If $\la$ is typical then 
there are no other $\mu \in\Z^n_+$ with $\wt(\mu) = \wt(\la)$.
So by the linkage principle \cite[Theorem 8.10]{BKI}, the induced module
$H^0(\la)$ is actually equal to $L(\la)$ in this case. 
Using this and the observation that $\u(\la)$ is an injective
$H$-supermodule for typical $\la$, 
we get from \cite[Theorem 7.5]{BKI} that $\Ext^1_G(L(\la),L(\mu)) 
= \Ext^1_G(L(\mu),L(\la)) = 0$
for all $\mu \in \Z^n_+$. 
In particular,  $L(\la)$ is injective,
hence $L(\la) = U(\la)$.
\end{proof}

Note in particular that the lemma shows that $U(\la)$ is self-dual
with respect to the duality $\tau$ in the case that $\la$ is typical.
Now suppose that $\la \in \Z^n_+$ is not typical.
Apply Procedure~\ref{alg2}
to get $\mu \in \Z^n_+$ and an operator $X_i \in \{E_i, F_i\}$.
Since this process reduces $\la$ to 
a typical weight in finitely many steps, we may assume inductively
that $U(\mu)$ has already been constructed, and that $U(\mu)
\cong U(\mu)^\tau$.
Just like in Lemma~\ref{props}, but applying Lemma~\ref{cor2} in place
of Lemma~\ref{cor}, we have that
$U(\la) \cong X_i U(\mu)$. Moreover,
$$
U(\la)^\tau \cong (X_i U(\mu))^\tau \cong X_i (U(\mu)^\tau) 
\cong X_i U(\mu) \cong U(\la),
$$
hence $U(\la)$ is also self-dual.
We obtain in this way an explicit
algorithm to construct all the injective indecomposables.
As a by-product we see that each $U(\la)$ is actually
self-dual with respect to the duality $\tau$, 
hence is isomorphic to the
{\em projective cover} of $L(\la)$.

Now we can prove the main result of the article. 
It shows that the map $\iota^*$ from (\ref{iotastardef}) maps
the canonical basis of $\F^n_{\Z}$ to the canonical basis of
$K(\mathcal C_n^*)$ given by the injective indecomposables,
and that the map $\iota$ from (\ref{iotadef}) maps the canonical basis
of $K(\mathcal C_n)$ given by the irreducible supermodules
to the dual canonical basis of $\E^n_{\Z}$.

\begin{Theorem}
For each $\la \in \Z^n_+$,
$\iota^* (U_\la(1)) = [U(\la)]$
and 
$\iota([L(\la)]) = L_\la(1)$.
\end{Theorem}

\begin{proof}
In view of Lemma~\ref{adj} and the facts that
 $[L(\la)]$ is dual to $[U(-w_0\la)]$ and
$L_\la(1)$ is dual to $U_{-w_0\la}(1)$,
it suffices to prove just that
$\iota^* (U_\la(1)) = [U(\la)]$.
If $\la$ is typical, 
then 
$U_\la(1) = F_\la(1)$
and $[U(\la)]  = [F(\la)]$, 
since there are no other $\mu \in \Z^n_+$ with $\wt(\mu) = \wt(\la)$.
So the result holds for typical weights.
The result follows in general because by definition
the map $\iota^*$ commutes with the operators $E_i, F_i$ and the algorithm
for constructing the $U(\la)$'s explained above exactly parallels
the algorithm for constructing the $U_\la(1)$'s from \ref{alg}.
\end{proof}

We get by the theorem and (\ref{ude}), (\ref{edef})
and (\ref{lfin}) respectively that 
\begin{align}
[E(\la)] &= \sum_{\mu \in \Z^n_+} u_{-w_0\la,-w_0\mu}(1) [L(\mu)],\label{edef2}\\
[L(\la)] &= \sum_{\mu \in \Z^n_+} l_{\mu,\la}(1) [E(\mu)],\label{lfin2}\\
[U(\la)] &= \sum_{\mu \in \Z^n_+} u_{\mu,\la}(1) [F(\mu)],
\end{align}
for each $\la \in \Z^n_+$.
Note finally that 
$u_{-w_0\la,-w_0\mu}(q) = u_{\la,\mu}(q)$, since 
$\omega(F_\mu) = F_{-w_0\mu}$ 
and $\omega(U_\mu) = U_{-w_0\mu}$ by Theorem~\ref{thmb}.
So comparing (\ref{edef2}) with (\ref{split})
gives that $d_{\mu,\la} = u_{\mu,\la}(1)$.
The Main Theorem stated in the introduction follows from
this statement and Theorem~\ref{mt}.

\Point{Conjectures}\label{mc}
To conclude the article, we make the following conjecture:
for each $\la \in \Z^n$, we have that
\begin{equation}\label{mc1}
[M(\la)] = \sum_{\mu \in \Z^n} t_{-\la,-\mu}(1) [L(\mu)].
\end{equation}
Comparing (\ref{tde}) and Corollary~\ref{fne},
this conjecture is equivalent to the statement that
\begin{equation}
[L(\la)] = \sum_{\mu \in \Z^n} l_{\mu,\la}(1) [M(\mu)],
\end{equation}
equality in $\widehat{K}(\mathcal O_n)$.
The $N_\la$'s in section 2 correspond to the supermodules
$N(\la)$ defined by
$N(\la) := U(\mathfrak g) \otimes_{U(\mathfrak b)} \tilde u(\la)$,
where $\tilde u(\la)$ denotes the projective cover of
$\u(\la)$ in the category of 
$\mathfrak h$-supermodules that are semisimple over $\mathfrak h_{\0}$.
Note at least that $[N(\la)]  = 2^{z(\la)} [M(\la)]$ in the
Grothendieck group, cf. (\ref{mdef}) and \cite[(7.1)]{BKI}.
The $T_\la$'s in section 2 should correspond to 
the {\em indecomposable tilting modules} $T(\la)$
in category $\mathcal O_n$, cf. \cite[Example 7.10]{Btilt}.
We recall that for $\la \in \Z^n$, $T(\la)$ is the supermodule
characterized uniquely up to isomorphism by the properties:
\begin{itemize}
\item[(1)] $T(\la) \in \mathcal O_n$ is indecomposable;
\item[(2)] $\Ext_{\mathcal O_n}^1(N(\mu),T(\la)) = 0$ for all $\mu \in \Z^n$;
\item[(3)] $T(\la)$ has a filtration where the subquotients are of the form 
$N(\mu)$ for
$\mu \in \Z^n$, starting with $N(\la)$ at the bottom.
\end{itemize}
Conjecture (\ref{mc1}) is equivalent
to the statement
\begin{equation}
[T(\la)] = \sum_{\mu \in \Z^n} t_{\mu,\la}(1) [N(\mu)],
\end{equation}
as follows by \cite[(7.12)]{Btilt}.


\begin{thebibliography}{BGG}
\parskip -0.5mm


\bibitem[Bg]{Berg}
G. Bergman,
The diamond lemma for ring theory,
{\em Advances Math.} {\bf 29} (1978), 178--218.

\bibitem[BGG]{BGG}
J. Bernstein, I. M. Gelfand and S. I. Gelfand,
A category of $\mathfrak g$-modules,
{\em Func. Anal. Appl.} {\bf 10} (1976), 87--92.



\bibitem[B1]{KL}
J. Brundan, 
Kazhdan-Lusztig polynomials and character formulae for the Lie superalgebra
$\mathfrak{gl}(m|n)$, preprint.

\bibitem[B2]{BII}
J. Brundan, Modular representations of the supergroup Q(n), II,
preprint.

\bibitem[B3]{Btilt}
J. Brundan, Tilting modules for Lie superalgebras,
preprint.



\bibitem[BK1]{BK}
J. Brundan and A. Kleshchev, Hecke-Clifford superalgebras,
crystals of type $A_{2\ell}^{(2)}$ and modular branching rules for 
$\widehat{S}_n$,
{\em Represent. Theory} {\bf 5} (2001), 317-403.

\bibitem[BK2]{BKI}
J. Brundan and A. Kleshchev, Modular representations of the supergroup Q(n), I,
to appear in {\em J. Algebra}.




























\bibitem[J1]{Jantzen}
J.~C. Jantzen, {\em Representations of algebraic groups}, Academic Press, 1986.

\bibitem[J2]{Jantzen2}
J.~C. Jantzen, {\em Lectures on quantum groups},
{Graduate Studies in Math.} {\bf 6}, Amer. Math. Soc., 1996.





\bibitem[JMO]{JMO}
N. Jing, K. Misra and M. Okado,
$q$-Wedge modules for quantized enveloping algebras of classical type,
{\em J. Algebra} {\bf 230} (2000), 518--539.

\bibitem[K1]{Kac}
V. G. Kac, Lie superalgebras, {\em Advances in Math.} {\bf 26} (1977), 8--96.

\bibitem[K2]{Kac2}
V. G. Kac, Characters of typical representations of classical Lie superalgebras,
{\em Commun. in Algebra} {\bf 5} (8) (1977), 889--897.

\bibitem[K3]{Kac3}
V. G. Kac, Representations of classical Lie superalgebras,
in: ``Differential geometrical methods in mathematical physics II'', Lecture Notes in Math. no. 676, pp. 597--626, Springer-Verlag, Berlin, 1978.


\bibitem[Ka1]{KaG}
M. Kashiwara, 
Global crystal bases of quantum groups,
{\em Duke Math. J.} {\bf 69} (1993), 455--485.

\bibitem[Ka2]{Ka}
M. Kashiwara,
On crystal bases,
{\em Proc. Canadian Math. Soc.} {\bf 16} (1995), 155--196.

\bibitem[KL]{KaL}
D. Kazhdan and G. Lusztig,
Representations of Coxeter groups and Hecke algebras, 
{\em Invent. Math.} {\bf 53} (1979), 165--184.






\bibitem[LLT]{LLT}
A. Lascoux, B. Leclerc and J.-Y. Thibon, Hecke algebras at roots of unity and crystal bases of quantum affine algebras, {\em Comm. Math. Phys.} {\bf 181} (1996), 205--263.

\bibitem[LT]{LT}
B. Leclerc and J.-Y. Thibon, $q$-Deformed Fock spaces and modular
representations of spin symmetric groups,
{\em J. Phys. A} {\bf 30} (1997), 6163--6176.


\bibitem[L1]{LuUnity}
G. Lusztig,
Quantum groups at roots of $1$,
{\em Geom. Ded.} {\bf 35} (1990), 89--114.

\bibitem[L2]{Lubook}
G. Lusztig,
{\em Introduction to quantum groups}, Progress in Math. {\bf 110}, Birkhauser, 
1993.

\bibitem[M]{Mac}
I. G. Macdonald, {\em Symmetric functions and Hall polynomials},
Oxford Mathematical Monographs, second edition, OUP, 1995.











\bibitem[P1]{penkov}
I. Penkov,
Characters of typical irreducible finite dimensional $\mathfrak q(n)$-supermodules,
{\em Func. Anal. Appl.} {\bf 20} (1986), 30--37.

\bibitem[P2]{p1}
I. Penkov,
Borel-Weil-Bott theory for classical Lie supergroups,
{\em Itogi Nauki i Tekhniki} {\bf 32} (1988),
71--124 (translation).

\bibitem[P3]{Penk}
I. Penkov, Generic representations of classical Lie superalgebras and
their localization, {\em Monatsh. Math.} {\bf 118} (1994), 267--313.

\bibitem[PS1]{PS}
I. Penkov and V. Serganova,
Cohomology of $G / P$ for classical complex Lie supergroups $G$
and characters of some atypical $G$-supermodules,
{\em Ann. Inst. Fourier} {\bf 39} (1989), 845--873.

\bibitem[PS2]{PS1} I. Penkov and V. Serganova, 
Characters of irreducible $G$-supermodules and cohomology of $G/P$ for the Lie supergroup $G=Q(N)$, {\em J. Math. Sci.} 
{\bf 84} (1997), no. 5, 1382--1412.

\bibitem[PS3]{PS2} I. Penkov and V. Serganova, Characters of finite dimensional 
irreducible $\mathfrak q(n)$-supermodules, {\em Lett. Math. Phys.} {\bf 40}
(1997), no. 2, 147--158.



\bibitem[S1]{sergeev}
A. N. Sergeev,
Tensor algebra of the identity representation as a module over the Lie
superalgebras $GL(n,m)$ and $Q(n)$, 
{\em Math. USSR Sbornik} {\bf 51} (1985), 419--427.

\bibitem[S2]{SergCent}
A. N. Sergeev,
The center of enveloping algebra for Lie
superalgebra $Q(n,{\mathbb C})$, 
{\em Lett. Math. Phys.} {\bf 7} (1983), 177--179.





\end{thebibliography}
\end{document}